\documentclass[sn-mathphys,Numbered]{sn-jnl}

\usepackage{graphicx}%
\usepackage{multirow}%
\usepackage{amsmath,amssymb,amsfonts}%
\usepackage{amsthm}%
\usepackage{mathrsfs}%
\usepackage[title]{appendix}%
\usepackage{xcolor}%
\usepackage{textcomp}%
\usepackage{manyfoot}%
\usepackage{booktabs}%
\usepackage{algorithm}%
\usepackage{algorithmicx}%
\usepackage[section]{placeins}
\usepackage{algpseudocode}%
\usepackage{listings}%
\usepackage{bm}%
\usepackage{subfigure}%
\theoremstyle{thmstyleone}%
\newtheorem{theorem}{Theorem}
\theoremstyle{thmstyletwo}%
\newtheorem{remark}{Remark}%
\newtheorem{lemma}{Lemma}
\theoremstyle{thmstylethree}%

\begin{document}

\title[DtN-AFEM for Acoustic-Elastic Interaction Problem]{An Adaptive Finite Element DtN Method for the Acoustic-Elastic Interaction Problem in Periodic Structures}


\author[1,2]{\fnm{Lei} \sur{Lin}}\email{linlei@amss.ac.cn}

\author*[1]{\fnm{Junliang} \sur{Lv}}\email{lvjl@jlu.edu.cn}
\equalcont{These authors contributed equally to this work.}


\affil*[1]{\orgdiv{School of Mathematics}, \orgname{Jilin University}, 
\orgaddress{\street{Qianjin Street}, \city{ Changchun}, \postcode{130012}, 
\state{Jilin Province}, \country{China}}}

\affil[2]{\orgdiv{ LSEC, Institute of Computational Mathematics and Scientific/Engineering Computing, Academy of Mathematics and Systems Science}, \orgname{Chinese Academy of Sciences},  \city{Beijing}, \country{China}}



\abstract{Consider a time-harmonic acoustic plane wave incident onto an elastic body with an unbounded periodic surface. The medium above the surface is supposed to be filled with a homogeneous compressible inviscid air/fluid of constant mass density, while the elastic body is assumed to be isotropic and linear. 
	By introducing the Dirichlet-to-Neumann (DtN) operators for acoustic and elastic waves simultaneously, the model is formulated as an acoustic-elastic interaction problem in periodic structures. Based on a duality argument, an a
	posteriori error estimate is derived for the associated truncated finite element approximation.
	The a posteriori error estimate consists of the finite element approximation error and the truncation error
	of two different DtN operators, where the latter decays exponentially with respect to the truncation parameter. Based on the a posteriori error, an adaptive finite element algorithm is proposed for solving the acoustic-elastic interaction problem in periodic structures. Numerical experiments are presented
	to demonstrate the effectiveness of the proposed algorithm.}




\maketitle

\backmatter

\section{Introduction}
Direct and inverse scattering form an important research area in science and engineering field\cite{Bao21Lin,Colton98,Li24,Chen08,Bao05,Lin25,Wang25}.
The scattering in periodic structures, which is known as diffraction gratings, has many significant applications in micro-optics including the design and fabrication of optical elements such as corrective lenses, anti-reflective interfaces, beam splitters, and sensors \cite{Bao13,Bao22,Niu24,Zhang18}. 
This paper is concerned with the scattering of a
time-harmonic acoustic plane wave by  an elastic body with an unbounded periodic surface.
The medium above the surface is supposed to be
filled with a homogeneous compressible inviscid air/fluid, and the
domain below the surface is occupied by an isotropic and linearly elastic solid body. Due to the external incident acoustic field, an elastic wave propagating downward is incited inside the solid, and the incident acoustic wave is scattered back into the air/fluid. Such a phenomenon can be modeled by the acoustic-elastic interaction problems in periodic structures, which have received wide attention due to their significant applications in
underwater acoustics, sonic, photonic crystals and ultrasonic nondestructive evaluation\cite{Claeys03,Declercq05}. A general mathematical framework for the acoustic-elastic interaction problem in both one- and two-dimensional periodic structures was established by variational methods in \cite{Hu16}.
Several computational approaches
have also been developed to numerically solve these problems, such as T matrix method \cite{Bishop93} and finite element method \cite{Hu16}.

Since the acoustic-elastic interaction problem in periodic structures is a class of scattering problems defined in an unbounded domain, it is necessary to truncate the unbounded domain into a bounded computational domain and impose an appropriate boundary condition on the artificial boundary for the application of the numerical methods, such as the finite difference method, the finite element method and so on. Such a boundary condition is named the transparent boundary condition (TBC), which is an important and active subject in the research area of wave propagation \cite{Bayliss80,Engquist77,Gachter03,Grote95,Grote04,Hagstrom99}. Another technique of  domain truncation is the perfectly
matched layer (PML), which was proposed by Berenger to solve the time-dependent Maxwell's equations \cite{Berenger94}. The main idea of the PML technique is using a layer of finite thickness fictitious medium
to surround
the computational domain such that
the waves coming from inside of the computational domain are absorbed. When the waves reach
the outer boundary of the PML region, their values are so small that the homogeneous
Dirichlet boundary conditions can be imposed on the boundary. Ever since Berenger's pioneering work,
the research on the PML has undergone a rapid development due to its effectiveness and simplicity. Various constructions of PML have been proposed
for different scattering problems, such as acoustic waves\cite{Chen05}, elastic waves \cite{Bramble10,Jiang18} and electromagnetic waves  \cite{Chew94,Hohage03}. 

The a posteriori error estimates are computable quantities measuring the errors between discrete problems and original problems, which are also essential in designing the adaptive finite element algorithm. Based on the a posteriori error estimates, the adaptive finite element methods
have the ability of error control and asymptotically optimal approximation property
\cite{Babuska78}. Therefore, the adaptive finite element methods have become a class of effective numerical methods for solving differential equations, especially for those problems whose solutions exist singularity or multiscale
phenomena. Combined with the PML technique, an efficient adaptive finite element
method was firstly developed in \cite{Chen03} for solving the diffraction grating
problem. It was shown that the a posteriori error estimate consists of the finite element discretization error and the PML
truncation error, where the latter decays exponentially with respect to the PML parameters, such
as the thickness of the layer and the medium properties. Due to its superior numerical behaviour, the adaptive finite element PML method was quickly developed to solve a variety of wave scattering problems, such as acoustic scattering\cite{Chen05,Jiang19}, elastic scattering\cite{Chen16,Jiang17-3,Jiang18}, electromagnetic scattering\cite{Bao10,Chen08,Zhou18} and acoustic-elastic interaction\cite{Jiang17-2}. 

As an alternative to the adaptive finite element PML method, the adaptive finite element DtN method (DtN-FEM)
has also been developed for solving
the obstacle scattering problems \cite{Bao21,Bao21DCDSB,Bao23,Jiang13,Jiang17,Lin24}, the diffraction grating problems \cite{Bao21-Jiang,Jiang22,Wang15,Yue23,Li20Yuan}, and the open cavity scattering problem \cite{Yuan20}. In this method, the transparent
boundary conditions are used to truncate the unbounded domain, and the layer
of artificial medium  used in PML is no longer needed. In fact, the transparent boundary conditions are given exactly. Therefore, the artificial boundary can
be chosen to closely surround the domain of interest, which can significantly reduce the size of the
computational domain. Since the transparent boundary conditions are usually given as infinite Fourier series, we need to truncate the infinite series into a sum of finite number of terms in practical computation. For theoretical analysis, the a posteriori error analysis technique of the adaptive finite element PML method
cannot be applied directly to the adaptive finite element DtN method. To overcome such an issue, a new duality argument was developed in \cite{Wang15} to obtain
the a posteriori error estimate between the solution of the original problem and that of
the discrete problem. Comparably, the a posteriori error estimates consist of
the finite element discretization error and the DtN truncation error, where the latter decays
exponentially with respect to the truncation parameter $N$. The numerical experiments also
demonstrate that the adaptive finite element DtN method is as effictive as the adaptive finite element PML method\cite{Jiang17,Jiang22,Wang15}. Recently, an adaptive finite volume element DtN method for the diffraction grating problem was also developed in
\cite{Wang22}.

In current paper, we present an adaptive DtN-FEM and carry out
its mathematical analysis for the acoustic-elastic interaction problem in periodic structures. The contribution is twofold: (1) give a
complete a posteriori error estimate; (2) develop an effective adaptive finite element
algorithm. This paper extends the work from the diffraction grating problem \cite{Wang15} to the acoustic-elastic interaction problem in periodic structures. Due to the complicated kinematic
and kinetic interface conditions imposed on the fluid-solid interface, the original model problem and the associated variational problem  of the acoustic-elastic interaction is much more complicated than those problems with single wave field. In the analysis of
the a posteriori error, we need to estimate two line integral terms defined on the fluid-solid interface, which
is different from the previous work about the single medium scattering. In addition, we have to consider the truncation errors of two different DtN operators simultaneously. We give the corresponding error
indicators of acoustic and elastic waves to reflect the distribution and magnitude
of the errors of different physical fields in the computation.

The outline of this paper is as follows. In Section \ref{model}, we introduce the model of the
acoustic-elastic interaction problem and its weak formulation by simultaneously introducing transparent boundaries of acoustic and elastic waves, respectively. The finite element discretization with two different truncated DtN operators is presented in Section \ref{fem}, while the a posteriori error estimate is derived by the duality
argument method in Section \ref{posterior}. In Section \ref{numerical}, we discuss the numerical implementation of our adaptive algorithm and
present some numerical experiments to demonstrate the effectiveness of the proposed
method. Eventually, the paper is concluded with some general remarks and directions for future work
in Section \ref{conclusion}.

For brevity, we use the notation $a \lesssim b$ to imply that $a \leq Cb$, where the positive constant
$C$ is independent of the truncation parameter of the DtN operators and
the mesh size of the triangulation.

\section{Model Problem}\label{model}
\begin{figure}
	\centering
	\includegraphics[width=0.5\textwidth]{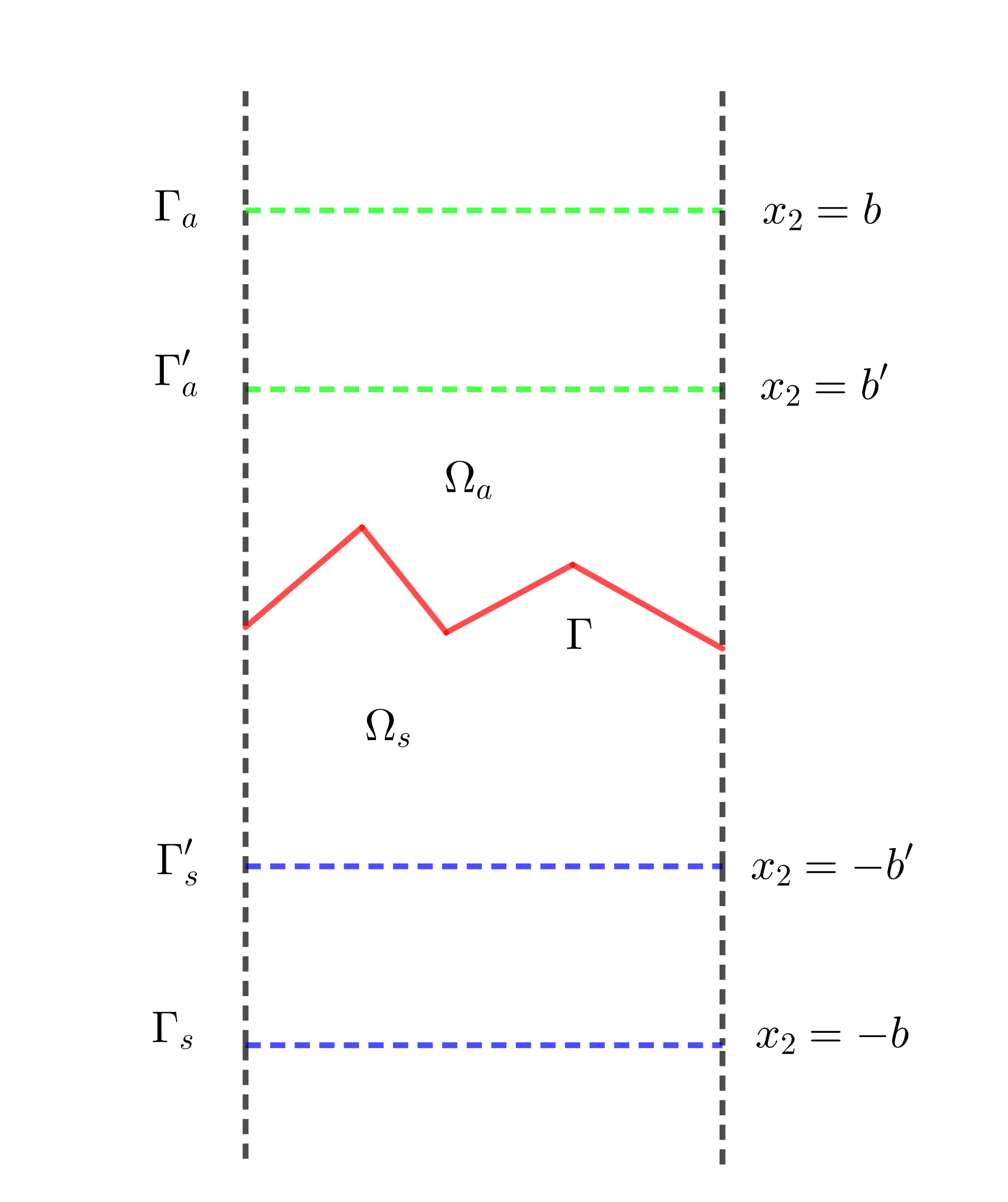}
	\caption{ A two-dimensional schematic of the problem geometry for the acoustic-elastic interaction in periodic structures}
	\label{model picture}
\end{figure}

Consider the acoustic wave incident onto the surface of an elastic solid, where the surface is assumed to
be periodic in the $x_1$-axis with the period $\Lambda$.  Due to the periodic structure of the elastic body, the problem
can be restricted into a single periodic cell, where 
$x_1 \in(0, \Lambda)$. Denote the surface of the elastic body by
\begin{equation*}
	\Gamma=\left\{\boldsymbol{x}=(x_1,x_2) \in \mathbb{R}^2: x_2=f(x_1), x_1 \in(0, \Lambda)\right\},
\end{equation*} 
where $f$
is a Lipschitz continuous function. The periodic
surface of the solid divides the whole space into the lower and upper domains. Denote
by $\Omega_+$ and $\Omega_-$ the upper and lower regions, respectively, where
\begin{align*}
	\Omega_+&=\{\boldsymbol{x} \in \mathbb{R}^2: f(x_1)<x_2, x_1 \in(0, \Lambda)\},\\
	\Omega_-&=\left\{\boldsymbol{x} \in \mathbb{R}^2: x_2<f(x_1), x_1 \in(0, \Lambda)\right\}.
\end{align*}
The domain $\Omega_+$ is filled with a homogeneous compressible inviscid fluid with a constant mass density $\rho_f > 0$. The domain $\Omega_-$ is occupied by an isotropic and linear elastic solid body determined through the Lam\'{e} constants $\lambda, \mu \in \mathbb{R}$ ($\mu>0$ and $\lambda+ \mu>0$) , and its mass density $\rho>0$. Let
\begin{align*}
	\Gamma_a&=\left\{\boldsymbol{x} \in \mathbb{R}^2: x_2=b, x_1 \in(0, \Lambda)\right\},\\
	\Gamma^{\prime}_a&=\left\{\boldsymbol{x} \in \mathbb{R}^2: x_2=b^{\prime}, x_1 \in(0, \Lambda)\right\},
\end{align*}
where $b$ and $b^{\prime}$ are constants satisfying $b>b^{\prime}>\max _{x_1 \in(0, \Lambda)} f(x_1)$. Similarly, we define
\begin{align*}
	\Gamma_s&=\left\{\boldsymbol{x} \in \mathbb{R}^2: x_2=-b, x_1 \in(0, \Lambda)\right\},\\
	\Gamma^{\prime}_s&=\left\{\boldsymbol{x} \in \mathbb{R}^2: x_2=-b^{\prime}, x_1 \in(0, \Lambda)\right\}.
\end{align*}
Let
\begin{align*}
	\Omega_a&=\left\{\boldsymbol{x} \in \mathbb{R}^2: f(x_1)<x_2<b, x_1 \in(0, \Lambda)\right\},\\
	\Omega_s&=\left\{\boldsymbol{x} \in \mathbb{R}^2: -b<x_2<f(x_1), x_1 \in(0, \Lambda)\right\}.
\end{align*}
The problem geometry is shown in Figure \ref{model picture}. 
Denote by $\bm{n}=(n_1,n_2)^{\top}$
the unit normal vector to $\Gamma$ directed into $\Omega_a$.

Let $p^i=e^{\mathrm{i}\left(\alpha x_1- \beta x_2 \right) }$ be the incident plane wave, where $\alpha=\kappa \sin \theta$, $\beta=\kappa \cos \theta$, and $
\theta \in \left(-\pi/2,\pi/2\right)$ is the incident angle.
Due to the wave interaction,
an elastic wave is induced inside the solid, and the scattered acoustic wave is generated in the
air/fluid. Such a phenomenon leads to an acoustic-elastic interaction problem. The incident plane wave $p^i$ satisfies the two-dimensional Helmholtz equation
\begin{equation*}
	\Delta p^i + \kappa^2 p^i = 0 \quad \text{in}~ \Omega_a, 
\end{equation*}
where $\kappa = \omega / c$ is the wavenumber, $\omega > 0$ is the angular frequency, and $c$ is the speed
of sound in the air/fluid. The total acoustic wave field $p^t$ also satisfies the Helmholtz equation
\begin{equation*}
	\Delta p^t + \kappa^2 p^t = 0 \quad \text{in}~ \Omega_a.
\end{equation*}
The total field $p^t$ consists of the incident field $p^i$ and the scattered field $p^s$, i.e.,
\begin{equation*}
	p^t=p^i+p^s \quad \text{in}~ \Omega_a. 
\end{equation*}
Let $H^1(\Omega_a)$  be the standard Sobolev space. Introduce the quasi-periodic functional space
\begin{align*}
	H_{\mathrm{qp}}^1(\Omega_a)&:=\left\{p \in H^1(\Omega_a): p(\Lambda, x_2)=p(0, x_2) e^{\mathrm{i} \alpha \Lambda}\right\}.
\end{align*}
Obviously, $H_{\mathrm{qp}}^1(\Omega_a)$ is the subspace of $H^1(\Omega_a)$.
For any function $p \in H_{\mathrm{qp}}^1(\Omega_a)$, it admits the Fourier expansion on $\Gamma_a$ as
\begin{align*}
	p(x_1, b)=\sum_{n \in \mathbb{Z}} p_n(b) e^{\mathrm{i} \alpha_n x_1}, \quad p_n(b)=\frac{1}{\Lambda} \int_0^{\Lambda} p(x_1, b) e^{-\mathrm{i} \alpha_n x_1} \mathrm{d} x_1,
\end{align*}
where $ \alpha_n=\alpha+n\left(2 \pi/\Lambda\right)$. Given $s \in \mathbb{R}$, the quasi-periodic trace functional space $H^s(\Gamma_a) $ is defined by
\begin{equation*}
	H^s(\Gamma_a)=\left\{p \in L^2(\Gamma_a):\|p\|_{H^s(\Gamma_a)}<\infty\right\},
\end{equation*}
where the trace norm $\|\cdot\|_{H^s(\Gamma_a)}$ is characterized by
\begin{equation*}
	\|p\|_{H^s(\Gamma_a)}=\left(\Lambda \sum_{n \in \mathbb{Z}}\left(1+\alpha_n^2\right)^s\left|p_n(b)\right|^2\right)^{1 / 2}.
\end{equation*}
It was shown in \cite{Wang15} that $p^s$ satisfies the upward Rayleigh expansion
\begin{equation}\label{Rayleigh1}
	p^s(x_1,x_2)=\sum_{n \in \mathbb{Z}} p^s_n(b) e^{\mathrm{i}\left(\alpha_n x_1+\beta_n(x_2-b)\right)}, \quad x_2>b,
\end{equation}
where
\begin{equation*}
	\beta_n= \begin{cases}\left(\kappa^2-\alpha_n^2\right)^{1 / 2}, & \left|\alpha_n\right|<\kappa, \\ \mathrm{i}\left(\alpha_n^2-\kappa^2\right)^{1 / 2}, & \left|\alpha_n\right|>\kappa .\end{cases}
\end{equation*}
To avoid the Wood's anomalies\cite{Bao22}, we assume that  $\kappa \neq\left|\alpha_n\right|$ . Taking the normal derivative of (\ref{Rayleigh1}) with respect to $x_2$ and then evaluating it at $x_2=b$, one can get
\begin{equation}\label{rayleigh}
	\partial_{x_2} p^s(x_1,b)=\sum_{n \in \mathbb{Z}} \mathrm{i} \beta_n p^s_n(b) e^{\mathrm{i} \alpha_n x_1}.
\end{equation}
Given a quasi-periodic function $p \in H_{\mathrm{qp}}^1(\Omega_a)$ which
admits a Fourier series expansion
\begin{align*}
	p(x_1, b)=\sum_{n \in \mathbb{Z}} p_n(b) e^{\mathrm{i} \alpha_n x_1},
\end{align*}
we define an acoustic DtN operator on $\Gamma_a$ as 
\begin{equation*}
	\mathscr{T}^a p(x_1,b) := \sum_{n \in \mathbb{Z}} \mathrm{i} \beta_n p_n(b) e^{\mathrm{i} \alpha_n x_1}.
\end{equation*}
By (\ref{rayleigh}), the following transparent
boundary condition can be imposed for $p^s$:
\begin{align}\label{DtN1}
	\partial_{x_2} p^s = \mathscr{T}^a p^s \quad \text{on}~ \Gamma_a.
\end{align}

The time-harmonic elastic transmitted wave is governed by the  Navier equation
\begin{equation}\label{navier equation}
	\Delta^{*} \bm{u} + \rho\omega^2 \bm{u} = 0 \quad \text{in}~ \Omega_s,
\end{equation}
where $\Delta^{*} := \mu \Delta + ( \lambda + \mu )\nabla \nabla\cdot$, and
$ \bm{u} = ( u_1, u_2 )^{\top} $ is the displacement of the elastic transmitted field. The traction operator $\bm{T}$ is defined by
\begin{equation*}
	\bm{Tu} := 2 \mu  \partial_{\bm{\nu}} \bm{u}  + \lambda \bm{\nu} \nabla \cdot \bm{u} + \mu
	\begin{bmatrix}
		\nu_2\left(\frac{\partial u_2}{\partial x_1} - \frac{\partial u_1}{\partial x_2}\right) \\
		\nu_1\left(\frac{\partial u_1}{\partial x_2} - \frac{\partial u_2}{\partial x_1}\right)
	\end{bmatrix},
\end{equation*}
where $\bm{\nu}=(\nu_1,\nu_2)^{\top}$ denotes the unit outward normal vector of $\Omega_s$.
To ensure the continuity of the normal component of the velocity, the kinematic interface  condition is assumed
\begin{equation}\label{transmission1}
	\partial_{\bm{n}} p^t  = \rho_f \omega^2 \bm{u} \cdot \bm{n} \quad \text{on}~ \Gamma.
\end{equation}
Besides, the following interface condition  is also assumed to ensure the continuity of the traction
\begin{equation}\label{transmission2}
	-p^t\bm{n} = \bm{Tu} \quad \text{on}~ \Gamma.
\end{equation}
For more details of the interface conditions (\ref{transmission1}) and (\ref{transmission2}), interested readers can refer to \cite{Hu16,Xu21}. 

Next we introduce the first generalized Betti's formula; see also \cite{Arens99}.
\begin{lemma}\label{Betti}
	Let $D \subset \mathbb{R}^2$ be a domain in which the divergence theorem holds. Then, for any
	vector fields $\bm{u} \in C^2(\overline{D})^2$ and $\bm{v}  \in C^1(\overline{D})^2$, the first generalized Betti formula holds:
	\begin{equation*}
		\int_{D} \Delta^*\bm{u} \cdot \bm{v}  \mathrm{d}\bm{x}
		+ \int_{D} \mathcal{E}_{\lambda,\mu} (\bm{u},\bm{v}) \mathrm{d}\bm{x} =  \int_{\partial D} \bm{Tu}\cdot \bm{v} \mathrm{d}s,
	\end{equation*}
	where
	\begin{align*}
		\mathcal{E}_{\lambda,\mu} (\bm{u},\bm{v})&=
		(2\mu+\lambda)\left(\frac{\partial u_1}{\partial x_1} \frac{\partial v_1}{\partial x_1}+
		\frac{\partial u_2}{\partial x_2} \frac{\partial v_2}{\partial x_2}
		\right)
		+\mu\left( \frac{\partial u_1}{\partial x_2} \frac{\partial v_1}{\partial x_2}+
		\frac{\partial u_2}{\partial x_1} \frac{\partial v_2}{\partial x_1} \right)\\
		&\quad+\lambda \left( \frac{\partial u_2}{\partial x_2} \frac{\partial v_1}{\partial x_1}+
		\frac{\partial u_1}{\partial x_1} \frac{\partial v_2}{\partial x_2} \right)
		+
		\mu \left( \frac{\partial u_2}{\partial x_1} \frac{\partial v_1}{\partial x_2}+
		\frac{\partial u_1}{\partial x_2} \frac{\partial v_2}{\partial x_1} \right) \\
		&=\lambda(\nabla \cdot \bm{u})(\nabla \cdot \bm{v}) 
		+\frac{\mu}{2}(\nabla \bm{u}+ \nabla \bm{u}^{\top} ):(\nabla \bm{v}+ \nabla \bm{v}^{\top}).
	\end{align*}
\end{lemma}
Let  $H^1(\Omega_s)$ be the standard Sobolev space and introduce its subspace
\begin{align*}
	H_{\mathrm{qp}}^1(\Omega_s)&:=\left\{u \in H^1(\Omega_s): u(\Lambda, x_2)=u(0, x_2) e^{\mathrm{i} \alpha \Lambda}\right\}.
\end{align*}
For any function $u \in H_{\mathrm{qp}}^1(\Omega_s)$, it admits the Fourier expansion on $\Gamma_s$ as
\begin{align*}
	u(x_1, -b)=\sum_{n \in \mathbb{Z}} u_n(-b) e^{\mathrm{i} \alpha_n x_1}, \quad u_n(-b)=\frac{1}{\Lambda} \int_0^{\Lambda} u(x_1,-b) e^{-\mathrm{i} \alpha_n x_1} \mathrm{d} x_1.
\end{align*}
The definition of the quasi-periodic trace functional space $H^s(\Gamma_s)$ and its trace norm $\|\cdot\|_{H^s(\Gamma_s)}$ can be given  similarly.
For any solution $\boldsymbol{u}$ of (\ref{navier equation}), it satisfies the Helmholtz decomposition
\begin{equation}\label{Helmholtz Decomposition}
	\boldsymbol{u}=\nabla \phi_1+\textbf{curl} \phi_2,
\end{equation}
where $\phi_j,j=1,2$ are scalar potential functions and $\textbf{curl} \phi_2=\left(\partial_{x_2} \phi_2,-\partial_{x_1} \phi_2\right)^{\top}$; see also \cite{Li22}. Substituting (\ref{Helmholtz Decomposition}) into (\ref{navier equation}), one can verify that $\phi_j, j=1,2$ satisfy the Helmholtz equation
\begin{equation}\label{Helmholtz}
	\Delta \phi_j+\kappa_j^2 \phi_j=0 \quad \text { in } \Omega_s,
\end{equation}
where $\kappa_1 = \omega\sqrt{\rho/(\lambda+2\mu)}$ and $\kappa_2=\omega\sqrt{\rho/\mu}$ are known as the compressional and shear wave numbers, respectively. Clearly, we have $\kappa_2>\kappa_1$.
Let $\phi_j$ be the solution of the Helmholtz equation (\ref{Helmholtz}) with the downward propogating wave condition. It is shown in \cite{Hu16} that $\phi_j$ is also a quasi-periodic function in the $x_1$ direction with period $\Lambda$ and satisfies the downward Rayleigh expansion
\begin{align}\label{Rayleigh2}
	\phi_j(x_1,x_2)=\sum_{n \in \mathbb{Z}} \phi_n^{(j)}(-b) e^{\mathrm{i}\left(\alpha_n x_1-\beta_n^{(j)}(x_2+b)\right)}, \quad x_2<-b.
\end{align}
where
\begin{equation}\label{beta_j}
	\beta_n^{(j)}= \begin{cases}\left(\kappa_j^2-\alpha_n^2\right)^{1 / 2}, & \left|\alpha_n\right|<\kappa_j, \\ \mathrm{i}\left(\alpha_n^2-\kappa_j^2\right)^{1 / 2}, & \left|\alpha_n\right|>\kappa_j .\end{cases}
\end{equation}
Similarly,
we assume that $\kappa_j \neq\left|\alpha_n\right|$ for $n \in \mathbb{Z}$.
Combining (\ref{Rayleigh2}) and the Helmholtz decomposition (\ref{Helmholtz Decomposition}) yields
\begin{equation}\label{u expression 1}
	\boldsymbol{u}(x_1,x_2)=\mathrm{i} \sum_{n \in \mathbb{Z}}\left[\begin{array}{c}
		\alpha_n \\
		-\beta_n^{(1)}
	\end{array}\right] \phi_n^{(1)}(-b) \mathrm{e}^{\mathrm{i}\left(\alpha_n x_1-\beta_n^{(1)}(x_2+b)\right)}+\left[\begin{array}{l}
		-\beta_n^{(2)} \\
		-\alpha_n
	\end{array}\right] \phi_n^{(2)}(-b) \mathrm{e}^{\mathrm{i}\left(\alpha_n x_1-\beta_n^{(2)}(x_2+b)\right)}.
\end{equation}
Besides, the displacement of the elastic wave $\boldsymbol{u}$ also satisfies the Fourier series expansion
\begin{align}\label{u expression 2}
	\boldsymbol{u}(x_1, -b)=\sum_{n \in \mathbb{Z}} \boldsymbol{u}_n(-b) \mathrm{e}^{\mathrm{i} \alpha_n x_1},
\end{align}
where
$\boldsymbol{u}_n(-b)=\left(u_n^{(1)}(-b),
u_n^{(2)}(-b)\right)^{\top}$.
From (\ref{u expression 1}) and (\ref{u expression 2}), we obtain a linear system of algebraic equations for $\phi_n^{(j)}(-b)$ :
\begin{equation}\label{algebraic equation}
	\left[\begin{array}{cc}
		\mathrm{i} \alpha_n & -\mathrm{i} \beta_n^{(2)} \\
		-\mathrm{i} \beta_n^{(1)} & -\mathrm{i} \alpha_n
	\end{array}\right]\left[\begin{array}{c}
		\phi_n^{(1)}(-b) \\
		\phi_n^{(2)}(-b)
	\end{array}\right]=\left[\begin{array}{c}
		u_n^{(1)}(-b) \\
		u_n^{(2)}(-b)
	\end{array}\right].
\end{equation}
Solving the above equations by Cramer's rule gives
\begin{align}
	& \phi_n^{(1)}(-b)=\frac{\mathrm{i}}{\chi_n}\left(-\alpha_n u_n^{(1)}(-b)+\beta_n^{(2)} u_n^{(2)}(-b)\right),\label{phi 1} \\
	& \phi_n^{(2)}(-b)=\frac{\mathrm{i}}{\chi_n}\left(\beta_n^{(1)} u_n^{(1)}(-b)+\alpha_n u_n^{(2)}(-b)\right),\label{phi 2}
\end{align}
where $\chi_n=\alpha_n^2+\beta_n^{(1)} \beta_n^{(2)}$.
Substituting (\ref{phi 1}) and (\ref{phi 2}) into (\ref{u expression 1}), we obtain the downward Rayleigh's expansion for $\bm{u}$: 
\begin{align}\label{rayleigh2}
	\boldsymbol{u}(x_1,x_2)&=  \sum_{n \in \mathbb{Z}} \frac{1}{\chi_n}\left[\begin{array}{cc}
		\alpha_n^2 & -\alpha_n \beta_n^{(2)} \\
		-\alpha_n \beta_n^{(1)} & \beta_n^{(1)} \beta_n^{(2)}
	\end{array}\right] \boldsymbol{u}_n(-b) \mathrm{e}^{\mathrm{i}\left(\alpha_n x_1-\beta_n^{(1)}(x_2+b)\right)} \nonumber\\
	& \qquad+\frac{1}{\chi_n}\left[\begin{array}{cc}
		\beta_n^{(1)} \beta_n^{(2)} & \alpha_n \beta_n^{(2)} \\
		\alpha_n \beta_n^{(1)} & \alpha_n^2
	\end{array}\right] \boldsymbol{u}_n(-b) \mathrm{e}^{\mathrm{i}\left(\alpha_n x_1-\beta_n^{(2)}(x_2+b)\right)},  \quad x_2<-b.
\end{align}
By the denifition of the traction operator and (\ref{rayleigh2}), we obtain the following transparent boundary condition for $\bm{u}$ on $\Gamma_s$ as
\begin{equation}\label{DtN2}
	\bm{Tu}=\mathscr{T}^s \boldsymbol{u}:=\sum_{n \in \mathbb{Z}} M_n\left(u_n^{(1)}(-b), u_n^{(2)}(-b)\right)^{\top} e^{\mathrm{i} \alpha_n x_1},
\end{equation}
where $\mathscr{T}^s$ is the elastic DtN operator and $M_n$ is a $2 \times 2$ matrix given by
\begin{equation}\label{M^n}
	M_n=\frac{\mathrm{i}}{\chi_n}\left[\begin{array}{cc}
		\omega^2 \beta_n^{(1)} & -2\mu \alpha_n \chi_n+\omega^2 \alpha_n \\
		2\mu \alpha_n \chi_n- \omega^2 \alpha_n & \omega^2 \beta_n^{(2)}
	\end{array}\right].
\end{equation}

By the DtN operators (\ref{DtN1}) and (\ref{DtN2}), the acoustic-elastic interaction problem in periodic structures can be formulated as the following coupled boundary value problem: Given $p^i$, find the quasi-periodic solutions $p^s$ and $\bm{u}$ such that
\begin{align}\label{DtNproblem}
	\Delta p^s + \kappa^2 p^s = 0 &\quad \text{in}~ \Omega_a, \nonumber\\
	\Delta^{*} \bm{u} + \rho\omega^2 \bm{u} = 0  &\quad \text{in}~ \Omega_s, \nonumber\\
	\partial_{\bm{n}} p^t  = \rho_f \omega^2 \bm{u} \cdot \bm{n} &\quad\text{on}~ \Gamma, \\
	-p^t\bm{n} =  \bm{Tu} &\quad \text{on}~ \Gamma,\nonumber\\
	\partial_y p^s = \mathscr{T}^a p^s&\quad \text{on}~  \Gamma_a, \nonumber\\
	\bm{Tu}  = \mathscr{T}^s \bm{u}& \quad 
	\text{on}~  \Gamma_s \nonumber. 
\end{align}
Let $\Omega = \Omega_a\cup\Omega_s$. Define 
\begin{equation*}
	\mathcal{H}^1(\Omega):=H^1_{\mathrm{qp}}(\Omega_a)\times H^1_{\mathrm{qp}}(\Omega_s)^2=\left\{U=(p, \boldsymbol{u}): p \in H^1_{\mathrm{qp}}\left(\Omega_a\right), \boldsymbol{u} \in H^1_{\mathrm{qp}}\left(\Omega_s\right)^2\right\}, 
\end{equation*}
which is endowed with the inner product
\begin{align*}
	(\bm{U},\bm{V})_{\mathcal{H}^1(\Omega)} :=&
	\int_{\Omega_s}\left(\lambda(\nabla \cdot \bm{u})(\nabla \cdot \overline{\bm{v}}) 
	+\frac{\mu}{2}(\nabla \bm{u}+ \nabla \bm{u}^{\top} ):(\nabla \overline{\bm{v}}+ \nabla \overline{\bm{v}}^{\top})+\bm{u} \cdot\overline{\bm{v}}
	\right)\text{d}\bm{x}\\
	&+\int_{\Omega_a} \left(\nabla p \cdot \nabla \overline{q} + p\overline{q} \right)\text{d}\bm{x}
\end{align*}
for any $\bm{U}=(p,\bm{u})$ and $\bm{V}=(q,\bm{v})$. Here, $A:B=\text{tr}(AB^{\top})$ is the Frobenius inner product of
square matrices $A$ and $B$, and $\nabla \bm{u}$ is the displacement gradient tensor, i.e.,
\begin{equation*}
	\nabla \bm{u}=
	\left[\begin{array}{cc}
		\frac{\partial u_1}{\partial x_1} & \frac{\partial u_1}{\partial x_2} \\
		\frac{\partial u_2}{\partial x_1} & \frac{\partial u_2}{\partial x_2}
	\end{array}\right].
\end{equation*}
One can prove that $\|\cdot\|_{\mathcal{H}^1(\Omega)}=\sqrt{(\cdot,\cdot)_{\mathcal{H}^1(\Omega)}}$ is a norm on $\mathcal{H}^1(\Omega)$. It can be easily verified that $\|\cdot\|_{\mathcal{H}^1(\Omega)}$ is equivalent to the standard $H^1$-norm; see the appendix in \cite{Xu21}.

It follows from (\ref{DtNproblem}), Green's formula and Lemma \ref{Betti} that
\begin{align*}
	0&=\int_{\Omega_a} \left( \Delta p^s \overline{q}+\kappa^2 p^s \overline{q}\right) \text{d}\bm{x} \\
	&=\int_{\Omega_a} \left(-\nabla p^s \cdot \nabla \overline{q} +\kappa^2 p^s\overline{q} \right)\text{d}\bm{x}-\int_{\Gamma} \rho_f\omega^2\bm{u}\cdot\bm{n} \overline{q} \text{d}s
	+\int_{\Gamma_a} \mathscr{T}^a p^s \overline{q}\text{d}s
	+\int_{\Gamma} \partial_{\bm{n}} p^i \overline{q} \text{d}s
\end{align*}
and
\begin{align*}
	0&=\int_{\Omega_s} \left(\Delta^*\bm{u} \cdot \bm{v}+ \rho\omega^2\bm{u} \cdot\overline{\bm{v}}
	\right)  \text{d}\bm{x}\\
	&=\int_{\Omega_s}\left(-\mathcal{E}_{\lambda,\mu} (\bm{u},\bm{v})+ \rho\omega^2\bm{u} \cdot\overline{\bm{v}}
	\right)\text{d}\bm{x}+\int_{\partial\Omega_s} \bm{Tu} \cdot \overline{\bm{v}} \text{d}s\\
	&=\int_{\Omega_s}\left(-\mathcal{E}_{\lambda,\mu} (\bm{u},\bm{v})+ \rho\omega^2\bm{u} \cdot\overline{\bm{v}}
	\right)\text{d}\bm{x}
	-\int_{\Gamma}p^s\bm{n} \cdot \overline{\bm{v}}\text{d}s
	-\int_{\Gamma} p^i\bm{n} \cdot \overline{\bm{v}}\text{d}s+\int_{\Gamma_s} \mathscr{T}^s \bm{u}\cdot \overline{\bm{v}}\text{d}s.
\end{align*}
Let $A: \mathcal{H}^1(\Omega)\times \mathcal{H}^1(\Omega) \to \mathbb{C}$ be the sesquilinear form as
\begin{equation}\label{defA}
	A(\bm{U},\bm{V})=a_1(\bm{u},\bm{v})+a_2(p^s,q)+a_3(\bm{u},q)+a_4(p^s,\bm{v})+b_1(p^s,q)+b_2(\bm{u},\bm{v}),
\end{equation}
where 
\begin{align*}
	a_1(\bm{u},\bm{v})&=\int_{\Omega_s}\left(\mathcal{E}_{\lambda,\mu} (\bm{u},\bm{v})- \rho\omega^2\bm{u} \cdot\overline{\bm{v}}
	\right)\text{d}\bm{x},\\
	a_2(p^s,q)&=\int_{\Omega_a} \left(\nabla p^s \cdot \nabla \overline{q} -\kappa^2 p^s\overline{q} \right)\text{d}\bm{x},\\
	a_3(\bm{u},q)&=\int_{\Gamma} \rho_f\omega^2\bm{u}\cdot\bm{n} \overline{q} \text{d}s,\\
	a_4(p^s,\bm{v})&=\int_{\Gamma}p^s\bm{n} \cdot \overline{\bm{v}}\text{d}s,\\
	b_1(p^s,q)&=-\int_{\Gamma_a} \mathscr{T}^a p^s \overline{q}\text{d}s,\\
	b_2(\bm{u},\bm{v})&=-\int_{\Gamma_s} \mathscr{T}^s \bm{u}\cdot \overline{\bm{v}}\text{d}s,
\end{align*}
and $\ell$ is a bounded linear functional dependent on $\left(\partial_{\bm{n}} p^i, p^i \right) \in H^{-1 / 2}(\Gamma) \times H^{1 / 2}(\Gamma) $, defined by
\begin{equation}\label{righthand}
	\ell (\bm{V})=\int_{\Gamma} \partial_{\bm{n}} p^i  \overline{q} - p^i\bm{n} \cdot \overline{\bm{v}}\text{d}s.
\end{equation}
Thus, we arrive at an equivalent variational formulation for (\ref{DtNproblem}): Given $p^i$, find  $\bm{U}=(p^s,\bm{u}) \in \mathcal{H}^1(\Omega)$ such that
\begin{equation}\label{varational}
	A(\bm{U}, \bm{V})=\ell(\bm{V})\quad ~\text{for all} ~ \bm{V}=(q,\bm{v}) \in \mathcal{H}^1(\Omega).
\end{equation}

Following the idea in reference \cite{Hsiao04}, the inf-sup stability (\ref{infsup}) can be obtained by the G\r{a}rding's inequality and the uniqueness result of (\ref{varational}). A similar strategy has also been used for the proof of
inf-sup condition of discrete sesquilinear form for acoustic-elastic interaction obstacle scattering problem in \cite{Xu21}.
Following \cite[Theorem 3]{Xu21}, we can prove that the sesquilinear form in (\ref{varational}) satisfies the G\r{a}rding's inequality. In addition, the uniqueness of the variational problem (\ref{varational}) has been proved in \cite{Hu16}.
Thus, there exists a constant $\varsigma$ such that the following inf-sup condition holds
\begin{align}\label{infsup}
	\sup_{0 \neq \bm{V} \in \mathcal{H}^1(\Omega)} \frac{\left|A(\bm{U},\bm{V}) \right|}{\|\bm{V}\|_{\mathcal{H}^1(\Omega)}} \ge \varsigma \|\bm{U}\|_{\mathcal{H}^1(\Omega)}
	\quad ~\text{for all} ~ \bm{U} \in \mathcal{H}^1(\Omega).
\end{align}
It follows from (\ref{righthand}), (\ref{varational}) and (\ref{infsup}) that
\begin{align}\label{regular1}
	\|\bm{U} \|_{\mathcal{H}^1(\Omega)} \lesssim
    \left\|p^i \right\|_{L^2(\Gamma)}+\left\|\partial_{\bm{n}} p^i \right\|_{L^2(\Gamma)}.
\end{align}

\section{Finite Element Approximation}\label{fem}
Denote by $\mathcal{M}_{h,s}$ and $\mathcal{M}_{h,a}$ the conforming triangulations of $\Omega_s$ and $\Omega_a$, respectively. Assume that 
$\mathcal{M}_{h} = \mathcal{M}_{h,s} \cup \mathcal{M}_{h,a}$.
Let $\hat{S}_h \subset H_{\mathrm{qp}}^1(\Omega_a)$ and $\overline{S}_h \subset H_{\mathrm{qp}}^1(\Omega_s)^2$ be the  conforming finite element spaces, i.e.,
\begin{align*}
	\hat{S}_h&:=\{ p_h \in C(\overline{\Omega}_a): p_h|_K \in P_m(K) ~\text{for all} ~K \in \mathcal{M}_{h,a}, p_h(\Lambda, x_2)=p_h(0, x_2) e^{\mathrm{i} \alpha \Lambda} \},\\
	\overline{S}_h&:=\{ \bm{u}_h \in C(\overline{\Omega}_s): \bm{u}_h|_K \in P_m(K)^2 ~\text{for all} ~K \in \mathcal{M}_{h,s},\bm{u}_h(\Lambda, x_2)=\bm{u}_h(0, x_2) e^{\mathrm{i} \alpha \Lambda} \},
\end{align*} 
where $m$ is a positive integer and $P_m(K)$ denotes the set of all polynomials of degree no more than $m$. Let $\mathcal{H}_h^1(\Omega) := \hat{S}_h \times \overline{S}_h \subset \mathcal{H}^1(\Omega)$.
The finite element approximation to the problem (\ref{varational}) reads as follows: Given $p^i$,  find $\bm{U}_h=(p_h^s,\bm{u}_h) \in \mathcal{H}_h^1(\Omega)$
such that
\begin{equation}\label{varationalFEM}
	A(\bm{U}_h, \bm{V}_h)=\ell(\bm{V}_h) \quad ~\text{for all} ~ \bm{V}_h=(q_h,\bm{v}_h) \in \mathcal{H}_h^1(\Omega),
\end{equation} 
where  
\begin{align*}
	A\left(\bm{U}_h,\bm{V}_h\right)&=a_1(\bm{u}_h,\bm{v}_h)+a_2(p_h^s,q_h)+a_3(\bm{u}_h,q_h)+a_4(p_h^s,\bm{v}_h)+b_1(p_h^s,q_h)+b_2(\bm{u}_h,\bm{v}_h),\\
	a_1(\bm{u}_h,\bm{v}_h)&=\int_{\Omega_s}\left(\mathcal{E}_{\lambda,\mu} (\bm{u}_h,\bm{v}_h)- \rho\omega^2\bm{u}_h \cdot\overline{\bm{v}}_h
	\right)\text{d}\bm{x},\\
	a_2(p_h^s,q_h)&=\int_{\Omega_a} \left(\nabla p_h^s \cdot \nabla \overline{q}_h -\kappa^2 p_h^s\overline{q}_h \right)\text{d}\bm{x},\\
	a_3(\bm{u}_h,q_h)&=\int_{\Gamma} \rho_f\omega^2\bm{u}_h\cdot\bm{n} \overline{q}_h \text{d}s,\\
	a_4(p_h^s,\bm{v}_h)&=\int_{\Gamma}p_h^s\bm{n} \cdot \overline{\bm{v}}_h\text{d}s,\\
	b_1(p_h^s,q_h)&=-\int_{\Gamma_a} \mathscr{T}^a p_h^s \overline{q}_h\text{d}s,\\
	b_2(\bm{u}_h,\bm{v}_h)&=-\int_{\Gamma_s} \mathscr{T}^s \bm{u}_h\cdot \overline{\bm{v}}_h\text{d}s,\\
	\ell (\bm{V}_h)&=\int_{\Gamma} \partial_{\bm{n}} p^i  \overline{q}_h - p^i\bm{n} \cdot \overline{\bm{v}}_h\text{d}s.
\end{align*}

In the above equations, the DtN operators $\mathscr{T}^a$ and $\mathscr{T}^s$ defined in (\ref{DtN1}) and (\ref{DtN2}) are given by infinite series.
In practical computation, it is necessary to truncate the operators $\mathscr{T}^a$ and $\mathscr{T}^s$ by taking finitely many terms of
the expansions so as to attain a feasible algorithm. Given sufficiently large constants $N_1,N_2$, we define the  
truncated DtN operators $\mathscr{T}^a_N$ and $\mathscr{T}^s_N$ as
\begin{align}
	\mathscr{T}^a_N p^s  &= \sum_{|n|\leq N} \mathrm{i} \beta_n p^s_n(b) e^{\mathrm{i} \alpha_n x_1} \quad \text{on}~ \Gamma_a,\label{trDtN1} \\ 
	\mathscr{T}^s_N \boldsymbol{u}&=\sum_{|n|\leq N} M_n\left(u_n^{(1)}(-b), u_n^{(2)}(-b)\right)^{\top} e^{\mathrm{i} \alpha_n x_1} \quad \text { on } \Gamma_s,\label{trDtN2}
\end{align}
where $N=\max(N_1,N_2)$.
Using the truncated DtN operators $\mathscr{T}^a$ and $\mathscr{T}^s$, one can get the truncated finite element approximation to the
problem (\ref{varational}): Given $p^i$, find $\bm{U}_h^N=(p_h^{s,N},\bm{u}_h^N) \in \mathcal{H}_h^1(\Omega)$
such that
\begin{equation}\label{varationalN}
	A^N(\bm{U}_h^N, \bm{V}_h)=\ell(\bm{V}_h) \quad ~\text{for all} ~ \bm{V}_h \in \mathcal{H}_h^1(\Omega),
\end{equation}
where 
\begin{align*}
	A^N\left(\bm{U}_h^N,\bm{V}_h\right)&=a_1(\bm{u}_h^N,\bm{v}_h)+a_2(p_h^{s,N},q_h)+a_3(\bm{u}_h^{N},q_h)+a_4(p_h^{s,N},\bm{v}_h)+b_1^N(p_h^{s,N},q_h)+b_2^N(\bm{u}_h^N,\bm{v}_h),\\
	a_1(\bm{u}_h^N,\bm{v}_h)&=\int_{\Omega_s}\left(\mathcal{E}_{\lambda,\mu} (\bm{u}_h^N,\bm{v}_h)- \rho\omega^2\bm{u}_h^N \cdot\overline{\bm{v}}_h
	\right)\text{d}\bm{x},\\
	a_2(p_h^{s,N},q_h)&=\int_{\Omega_a} \left(\nabla p_h^{s,N} \cdot \nabla \overline{q}_h -\kappa^2 p_h^{s,N}\overline{q}_h \right)\text{d}\bm{x},\\
	a_3(\bm{u}_h^{N},q_h)&=\int_{\Gamma} \rho_f\omega^2\bm{u}_h^{N}\cdot\bm{n} \overline{q}_h \text{d}s,\\
	a_4(p_h^{s,N},\bm{v}_h)&=\int_{\Gamma}p_h^{s,N}\bm{n} \cdot \overline{\bm{v}}_h\text{d}s,\\
	b_1^N(p_h^{s,N},q_h)&=-\int_{\Gamma_a} \mathscr{T}_N^a p_h^{s,N} \overline{q}_h\text{d}s,\\
	b_2(\bm{u}_h^N,\bm{v}_h)&=-\int_{\Gamma_s} \mathscr{T}_N^s \bm{u}_h^{N} \overline{\bm{v}}_h\text{d}s.
\end{align*}

We refer readers to \cite{Hu16} for the well-posedness and the a priori error analysis of the problem (\ref{varationalN}).
In this work, we mainly focus on the a posteriori error estimate and the associated adaptive algorithm. Thus, we assume that the discrete problem 
(\ref{varationalN}) has a unique solution $\bm{U}_h^N \in \mathcal{H}_h^1(\Omega)$.

\section{The A Posteriori Error Analysis}\label{posterior}
For any $K \in \mathcal{M}_{h}$, denote by $h_{K}$ its diameter. Let $\mathcal{B}_h$ denote the set of all the edges of $K$. For
any $e \in \mathcal{B}_h$, denote by $h_e$ its length. For $K \in \mathcal{M}_{h,a}$ and $K \in \mathcal{M}_{h,s}$, denote the jump residual across the edge of $K$ by $J_{e,a}$ and $\bm{J_{e,s}}$, respectively.

For any interior edge $e$ which is the common side of $K_1, K_2\in  \mathcal{M}_{h,a}$, define the jump residual across $e$ as 
\begin{equation*}
	J_{e,a}=-\left(\nabla p_h^{s,N}|_{K_1} \cdot \bm{\nu}_1+\nabla p_h^{s,N}|_{K_2} \cdot \bm{\nu}_2 \right),
\end{equation*}
where $\bm{\nu}_j$ is the unit outward normal vector on the boundary of $K_j,j=1,2$. For any
boundary edge $e \in \Gamma_a$, define the jump residual
\begin{equation*}
	J_{e,a}=2\left(\mathscr{T}^a_N p_h^{s,N}-\nabla p_h^{s,N} \cdot \bm{\nu}\right),
\end{equation*}
where $\bm{\nu}$ is the unit outward normal on $ \Gamma_a$.
For any
boundary edge $e \in \Gamma$, define the jump residual
\begin{align*}
	J_{e,a}&=2\left( \partial_{\bm{n}} (p^i+p_h^{s,N}) -\rho_f\omega^2\bm{u}_h^{N} \cdot \bm{n}
	\right),\\
	\bm{J_{e,s}}&=-2\left( (p^i+p_h^{s,N})\bm{n} + \bm{Tu}_h^N
	\right).
\end{align*}
For any interior edge $e$ which is the common side of $K_1, K_2\in  \mathcal{M}_{h,s}$, define the jump residual across $e$ as 
\begin{equation*}
	\bm{J_{e,s}}=-\left(\bm{Tu}_h^N|_{K_1}+\bm{Tu}_h^N|_{K_2} \right).
\end{equation*}
For any
boundary edge $e \in \Gamma_s$, define the jump residual
\begin{equation*}
	\bm{J_{e,s}}=2\left(\mathscr{T}^s_N \bm{u}_h^{N}-\bm{Tu}_h^N\right).
\end{equation*}
For any boundary edge on the left line segment of $\partial \Omega$, i.e., $e \in\{x_1=0\} \cap \partial K_1$ for some $K_1 \in \mathcal{M}_h$, and its corresponding edge on the right line segment of $\partial \Omega$, i.e., $e^{\prime} \in\{x_1=\Lambda\} \cap \partial K_2$ for some $K_2 \in \mathcal{M}_h$, define the jump residual as

\begin{align*}
	J_{e,a}&=-\left(\nabla p_h^{s,N}|_{K_1} \cdot \bm{\nu}_1+e^{-\mathrm{i} \alpha \Lambda}\nabla p_h^{s,N}|_{K_2} \cdot \bm{\nu}_2 \right),\\
	J_{e^{\prime},a}&=-\left(e^{\mathrm{i} \alpha \Lambda}\nabla p_h^{s,N}|_{K_1} \cdot \bm{\nu}_1+\nabla p_h^{s,N}|_{K_2} \cdot \bm{\nu}_2 \right),\\
	\bm{J_{e,s}} & =-\left(\bm{Tu}_h^N|_{K_1}+e^{-\mathrm{i} \alpha \Lambda}\bm{Tu}_h^N|_{K_2}\right),\\
	\bm{J_{e^{\prime},s}} & =-\left(e^{\mathrm{i} \alpha \Lambda}\bm{Tu}_h^N|_{K_1}+\bm{Tu}_h^N|_{K_2}\right).
\end{align*}

For $K \in \mathcal{M}_{h,a}$ and $K \in \mathcal{M}_{h,s}$, denote the local error estimators by
$\eta_{K,a}$ and $\eta_{K,s}$, which are defined respectively by
\begin{equation}\label{estimator1}
	\eta_{K,a}=h_K\left\| \mathcal{R}_a p_h^{s,N}\right\|_{L^2(K)}+
	\left( \frac{1}{2}\displaystyle \sum_{e \in \partial K} h_e\| J_{e,a}\|_{L^2(e)}^2
	\right)^{1/2},
\end{equation}
\begin{equation}\label{estimator2}
	\eta_{K,s}=h_K\left\| \mathcal{R}_s \bm{u}_h^{N}\right\|_{L^2(K)}+
	\left( \frac{1}{2}\displaystyle \sum_{e \in \partial K} h_e\| \bm{J_{e,s}}\|_{L^2(e)}^2
	\right)^{1/2},
\end{equation}
where  $\mathcal{R}_{a}:=\Delta+\kappa^2$
and $\mathcal{R}_{s}:=\Delta^*+\rho\omega^2$.

Denote the error by $\bm{\Psi}=(\xi,\bm{\zeta}) \in \mathcal{H}^1(\Omega)$, where $\xi:= p^s-p_h^{s,N}, \bm{\zeta}:= \bm{u}-\bm{u}_h^{N}$. Introduce a dual problem to the original scattering problem:
Find $\bm{\Phi}=(\gamma,\bm{w}) \in \mathcal{H}^1(\Omega)$ such that
\begin{equation}\label{dual}
	A(\bm{V},\bm{\Phi})=(\bm{V},\bm{\Psi})\quad \text{for all}~ \bm{V}=(q,\bm{v}) \in \mathcal{H}^1(\Omega),
\end{equation}
where
\begin{align*}
	A(\bm{V},\bm{\Phi})=A(q,\bm{v};\gamma,\bm{w}), 
\end{align*}
and
\begin{equation*}
	(\bm{V},\bm{\Psi})=(q,\bm{v};\xi,\bm{\zeta}).
\end{equation*}
It can be verified that $\bm{\Phi}$ is the weak solution to the following boundary value problem
\begin{align*}
	\Delta \gamma + \kappa^2 \gamma = -\xi &\quad \text{in}~ \Omega_a, \\
	\Delta^{*} \bm{w} +  \rho\omega^2 \bm{w} = -\bm{\zeta} &\quad  \text{in} ~\Omega_s, \\
	\partial_{\bm{n}} \gamma  = \bm{w} \cdot \bm{n} &\quad \text{on}~ \Gamma, \\
	-\rho_f \omega^2 \gamma\bm{n} = \bm{Tw} &\quad  \text{on}~ \Gamma, \\
	\partial_{y} \gamma = \mathscr{T}^{a,*} \gamma &\quad \text{on}~ \Gamma_a,\\
	\bm{Tw}  = \mathscr{T}^{s,*} \bm{w}& \quad 
	\text{on}~  \Gamma_s, 
\end{align*}
where the adjoint DtN operators $\mathscr{T}^{a,*}$ and $\mathscr{T}^{s,*}$ are defined by
\begin{align*}
	\mathscr{T}^{a,*}\gamma = &-\sum_{n \in \mathbb{Z}} \mathrm{i} \overline{\beta_n} \gamma_n(b) e^{\textup{i}\alpha_n x_1},\\
	\mathscr{T}^{s,*} \bm{w}= & =\sum_{n \in \mathbb{Z}}-\frac{\mathrm{i}}{\overline{\chi}_n}\left[\begin{array}{cc}
		\omega^2 \overline{\beta_n^{(1)}} & -\omega^2 \alpha_n+2\mu \alpha_n \overline{\chi}_n \\
		-2\mu \alpha_n \overline{\chi}_n+\omega^2 \alpha_n & \omega^2 \overline{\beta_n^{(2)}}
	\end{array}\right] \left(w_n^{(1)}(-b), w_n^{(2)}(-b)\right)^{\top} e^{\mathrm{i} \alpha_n x_1},
\end{align*}
respectively.
Like the original problem (\ref{varational}), we may show that
the dual problem (\ref{dual}) has a unique weak solution, which satisfies 
\begin{equation}\label{regular2}
	\|\bm{\Phi}\|_{\mathcal{H}^1(\Omega)} \lesssim \|\bm{\Psi} \|_{L^2(\Omega)}.
\end{equation}

The following lemma gives several energy representations of the error 
$\bm{\Psi}$.
\begin{lemma}
	Let $(p^s,\bm{u}), (p_h^{s,N},\bm{u}_h^N)$ and $\bm{\Phi}=(\gamma,\bm{w})$ be the solutions to the problems (\ref{varational}),  (\ref{varationalN}) and (\ref{dual}), respectively.
	We have
	\begin{align}\label{L2}
		\|\bm{\Psi}\|_{L^2(\Omega)}^2=&A(\bm{\Psi},\bm{\Phi})+\langle(\mathscr{T}^a-\mathscr{T}^a_N)\xi,\gamma\rangle_{\Gamma_a}+\langle(\mathscr{T}^s-\mathscr{T}^s_N)\bm{\zeta},\bm{w}\rangle_{\Gamma_s}\nonumber\\
		&-\langle(\mathscr{T}^a-\mathscr{T}^a_N)\xi,\gamma\rangle_{\Gamma_a}-\langle(\mathscr{T}^s-\mathscr{T}^s_N)\bm{\zeta},\bm{w}\rangle_{\Gamma_s},
	\end{align}
	\begin{align}
		\|\bm{\Psi}\|_{\mathcal{H}^1(\Omega)}^2=&\Re\left(A( \bm{\Psi},\bm{\Psi} ) + \langle(\mathscr{T}^a-\mathscr{T}^a_N)\xi,\xi\rangle_{\Gamma_a} + \langle(\mathscr{T}^s-\mathscr{T}^s_N)\bm{\zeta},\bm{\zeta}\rangle_{\Gamma_s}
		\right)\nonumber\\
		&+\Re\langle\mathscr{T}^a_N\xi,\xi\rangle_{\Gamma_a}+\Re\langle\mathscr{T}^s_N\bm{\zeta},\bm{\zeta}\rangle_{\Gamma_s}
		+(\rho\omega^2+1)\|\bm{\zeta} \|_{L^2(\Omega_s)}^2 \nonumber\\
		&+ (\kappa^2+1) \|\xi \|_{L^2(\Omega_a)}^2-\Re\left(a_3(\bm{\zeta},\xi)+a_4(\xi,\bm{\zeta})\right),\label{H1}
	\end{align}
	\begin{align}\label{error}
		&\quad A(\bm{\Psi},\bm{V})+\langle(\mathscr{T}^a-\mathscr{T}^a_N)\xi,q\rangle_{\Gamma_a}+\langle(\mathscr{T}^s-\mathscr{T}^s_N)\bm{\zeta},\bm{v}\rangle_{\Gamma_s}\nonumber\\
		&=\ell(q-q_h,\bm{v}-\bm{v}_h)-A^N\left(p_h^{s,N},\bm{u}_h^N;q-q_h,\bm{v}-\bm{v}_h\right)\nonumber\\
		&\quad+\langle(\mathscr{T}^a-\mathscr{T}^a_N)p^s,q\rangle_{\Gamma_a}+\langle(\mathscr{T}^s-\mathscr{T}^s_N)\bm{u},\bm{v}\rangle_{\Gamma_s},\quad \forall \boldsymbol{V} \in \mathcal{H}^1(\Omega), \boldsymbol{V}_h \in \mathcal{H}_h^1(\Omega).
	\end{align}
\end{lemma}
\begin{proof}
	Taking $\bm{V} =\bm{\Psi} $ in (\ref{dual})
	gives (\ref{L2}).
	It follows from (\ref{defA}) that
	\begin{equation*}
		A(\bm{\Psi},\bm{\Psi})=a_1(\bm{\zeta},\bm{\zeta})+a_2(\xi,\xi)+a_3(\bm{\zeta},\xi)+a_4(\xi,\bm{\zeta})+b_1(\xi,\xi)+b_2(\bm{\zeta},\bm{\zeta}),
	\end{equation*}
	which gives
	\begin{align*}
		\|\bm{\Psi}\|_{\mathcal{H}^1(\Omega)}^2=&A( \bm{\Psi},\bm{\Psi} ) +\langle\mathscr{T}^a\xi,\xi\rangle_{\Gamma_a}
		+\langle\mathscr{T}^s\bm{\zeta},\bm{\zeta}\rangle_{\Gamma_s}+(\rho\omega^2+1)\|\bm{\zeta} \|_{L^2(\Omega_s)}^2\\ 
		&+ (\kappa^2+1) \|\xi \|_{L^2(\Omega_a)}^2
		-a_3(\bm{\zeta},\xi)-a_4(\xi,\bm{\zeta}).
	\end{align*}
	Taking the real parts on both sides of the above equation yields (\ref{H1}).
	It follows from (\ref{varational}) and (\ref{varationalN}) that
	\begin{align*}
		A(\bm{\Psi},\bm{V})&=A\left( p^s-p_h^{s,N},\bm{u}-\bm{u}_h^N ; q-q_h,\bm{v}-\bm{v}_h\right)
		+A\left( p^s-p_h^{s,N},\bm{u}-\bm{u}_h^N ; q_h,\bm{v}_h\right)\\
		&=A\left( p^s,\bm{u}; q-q_h,\bm{v}-\bm{v}_h\right)
		-A\left( p_h^{s,N},\bm{u}_h^N;q-q_h,\bm{v}-\bm{v}_h\right)\\
		&\quad+A\left( p^s-p_h^{s,N},\bm{u}-\bm{u}_h^N;q_h,\bm{v}_h\right)\\
		&=\ell(q-q_h,\bm{v}-\bm{v}_h)-A^N\left(p_h^{s,N},\bm{u}_h^N ;q-q_h,\bm{v}-\bm{v}_h \right)\\
		&\quad+A^N\left(p_h^{s,N},\bm{u}_h^N ;q-q_h,\bm{v}-\bm{v}_h \right)
		-A\left( p_h^{s,N},\bm{u}_h^N;q-q_h,\bm{v}-\bm{v}_h\right)\\
		&\quad +A\left( p^s,\bm{u};q_h,\bm{v}_h\right)
		-A \left(p_h^{s,N},\bm{u}_h^N;q_h,\bm{v}_h\right).
	\end{align*}
	Since 	
	$A\left( p^s,\bm{u};q_h,\bm{v}_h\right)=\ell(q_h,\bm{v}_h)=A^N\left(p_h^{s,N},\bm{u}_h^N ;q_h,\bm{v}_h \right)$,
	we have
	\begin{align*}
		A(\bm{\Psi},\bm{V})&=\ell(q-q_h,\bm{v}-\bm{v}_h)-A^N\left(p_h^{s,N},\bm{u}_h^N;q-q_h,\bm{v}-\bm{v}_h \right)\\
		&\quad+A^N\left(p_h^{s,N},\bm{u}_h^N ;q,\bm{v} \right)
		-A\left(p_h^{s,N},\bm{u}_h^N;q,\bm{v}\right)\\
		&=\ell(q-q_h,\bm{v}-\bm{v}_h)-A^N\left(p_h^{s,N},\bm{u}_h^N;q-q_h,\bm{v}-\bm{v}_h \right)\\
		&\quad +\langle\left(\mathscr{T}^a-\mathscr{T}^a_N \right) p_h^{s,N}, q\rangle_{\Gamma_a}+\langle\left(\mathscr{T}^s-\mathscr{T}^s_N \right) \bm{u}_h^N, \bm{v}\rangle_{\Gamma_s}\\
		&=\ell(q-q_h,\bm{v}-\bm{v}_h)-A^N\left(p_h^{s,N},\bm{u}_h^N;q-q_h,\bm{v}-\bm{v}_h \right)\\
		&\quad-\langle\left(\mathscr{T}^a-\mathscr{T}^a_N \right) \xi, q\rangle_{\Gamma_a}+\langle\left(\mathscr{T}^a-\mathscr{T}^a_N \right) p^s, q\rangle_{\Gamma_a}\\
		&\quad-\langle\left(\mathscr{T}^s-\mathscr{T}^s_N \right) \bm{\zeta}, \bm{v}\rangle_{\Gamma_s}+\langle\left(\mathscr{T}^s-\mathscr{T}^s_N \right) \bm{u}, \bm{v}\rangle_{\Gamma_s},
	\end{align*}
	which implies (\ref{error}).
\end{proof}

Here, we present some technical estimates, which are used in our analysis of the a posterior error. For the sake of brevity, we have omitted the proof and provided the related literature for interested readers.
\begin{lemma}[See \cite{Jiang17-3}]\label{chi_n bound }
	For any $ n \in \mathbb{Z}$, we have  $\kappa_1^2<\left|\chi_n\right|<\kappa_2^2$.
\end{lemma}
\begin{remark}\label{M^n bound}
	For sufficiently large $|n|$, it can verified that $\alpha_n$ has an order of $n$.
	It follows from Lemma \ref{chi_n bound } that for sufficiently large $|n|$, we have
	\begin{equation*}
		|M_n^{(ij)}| \lesssim |n|,\quad i,j=1,2.
	\end{equation*}
\end{remark}
The following lemma is concerned with the trace of the solution to (\ref{varational}).
\begin{lemma}[\text{See} \cite{Chen03}]\label{trace}
	For any $p \in H_{\mathrm{qp}}^1(\Omega_a)$ and $u \in H_{\mathrm{qp}}^1(\Omega_s)$, there holds
	\begin{align*}
		&\| p\|_{H^{1/2}(\Gamma_a)} \lesssim \| p\|_{H^1(\Omega_a)}, \quad
		\| u\|_{H^{1/2}(\Gamma_s)} \lesssim \| u\|_{H^1(\Omega_s)},\\
		&\| p\|_{H^{1/2}(\Gamma_a^{\prime})} \lesssim \| p\|_{H^1(\Omega_a)}, \quad
		\| u\|_{H^{1/2}(\Gamma_s^{\prime})} \lesssim \| u\|_{H^1(\Omega_s)}.
	\end{align*}
\end{lemma}

The following Lemma \ref{hatp} is crucial in deriving the truncation error.
\begin{lemma}[\text{See} \cite{Wang15}]\label{hatp}
	Let $\bm{U}=(p^s,\bm{u})$ be the solution to (\ref{varational}). For sufficiently large $|n|$, we have
	\begin{equation*}
		\left|p^s_n(b)\right| \lesssim e^{-\left(b-b^{\prime}\right)\left|\beta_n\right|}  \left|p^s_n(b^{\prime})\right|.
	\end{equation*}
\end{lemma}

Besides, we need the following lemma to estimate two line integral terms defined on the fluid-solid interface.
\begin{lemma}[\text{See}\cite{Brenner08}]\label{trace theorem}
	Suppose that $D$ has a Lipschitz boundary, and that $p$ is
	a real number in the range $1 \leq p < \infty$. Then,  there holds
	\begin{equation*}
		\|v\|_{L^p(\partial D)} \lesssim \|v\|_{L^p(D)}^{1-1/p}\|v\|_{W_p^1(D)}^{1/p},
		\quad  \forall v \in W_p^1(D).
	\end{equation*}
\end{lemma}
The following lemma gives an upper bound estimate of (\ref{error}).
\begin{lemma}\label{estimate}
	For any $\bm{V} \in \mathcal{H}^1(\Omega)$, we have 
	\begin{align*}
		&\quad \left| A(\bm{\Psi},\bm{V})+\langle(\mathscr{T}^a-\mathscr{T}^a_N)\xi,q\rangle_{\Gamma_a}+\langle\left(\mathscr{T}^s-\mathscr{T}^s_N \right) \bm{\zeta}, \bm{v}\rangle_{\Gamma_s}\right|\\
		&\lesssim\left( \left( \displaystyle \sum_{K \in \mathcal{M}_{h,a}} \eta_{K,a}^2+
		\displaystyle \sum_{K \in \mathcal{M}_{h,s}} \eta_{K,s}^2 \right)^{1/2}+ \Theta \left(\left\|p^i \right\|_{L^2(\Gamma)}+\left\|\partial_{\bm{n}} p^i \right\|_{L^2(\Gamma)}\right)
		\right)\| \bm{V}\|_{\mathcal{H}^1(\Omega)},
	\end{align*}
	where 
	\begin{equation*}
		\Theta = \max\left(\max_{|n| > N}\left(e^{-\left(b-b^{\prime}\right)\left|\beta_n\right|}\right),\max _{|n|>N}\left(|n| e^{- \left(b-b^{\prime}\right)\left|\beta_n^{(2)}\right|}\right) \right).
	\end{equation*}
\end{lemma}
\begin{proof}
	Let
	\begin{align*}
		&J_1=\ell(q-q_h,\bm{v}-\bm{v}_h)-A^N\left(p_h^{s,N},\bm{u}_h^N;q-q_h,\bm{v}-\bm{v}_h\right),\nonumber\\
		&J_2=\langle(\mathscr{T}^a-\mathscr{T}^a_N)p^s,q\rangle_{\Gamma_a},\\
		&J_3=\langle(\mathscr{T}^s-\mathscr{T}^s_N)\bm{u},\bm{v}\rangle_{\Gamma_s},
	\end{align*}
	where $(q_h,\bm{v}_h) \in \mathcal{H}_h^1(\Omega)$. It follows from (\ref{error}) that
	\begin{equation*}
		A(\bm{\Psi},\bm{V})+\langle(\mathscr{T}^a-\mathscr{T}^a_N)\xi,q\rangle_{\Gamma_a}+\langle\left(\mathscr{T}^s-\mathscr{T}^s_N \right) \bm{\zeta}, \bm{v}\rangle_{\Gamma_s}=J_1+J_2+J_3.
	\end{equation*}
	By the definition of the sesquilinear form (\ref{varationalN}), we have
	\begin{align*}
		J_1=&\displaystyle \sum_{K \in \mathcal{M}_{h,a}}\left( \int_{K} 
		\left(-\nabla p_h^{s,N} \cdot \nabla(\overline{q}-\overline{q}_h) +\kappa^2 p_h^{s,N}(\overline{q}-\overline{q}_h)\right) \text{d}\bm{x}
		\right.\\
		&\left.\qquad\qquad + \displaystyle \sum_{e \in \partial_K \cap \Gamma}
		\int_{e}\left(  \partial_{\bm{n}} p^i 
		-\rho_f \omega^2 \bm{u}_h^N \cdot \bm{n}
		\right)(\overline{q}-\overline{q}_h)\text{d}s
		+\sum_{e \in \partial_K \cap \Gamma_a} \int_{e} \mathscr{T}^a_N p_h^{s,N}(\overline{q}-\overline{q}_h)\text{d}s \right)\\
		&+\sum_{K \in \mathcal{M}_{h,s}}\left(\int_{K}
		-\mathcal{E}_{\lambda,\mu}\left( \bm{u}_h^N,\bm{v}-\bm{v}_h\right)+ \rho\omega^2\bm{u}_h^N\cdot(\overline{\bm{v}}-\overline{\bm{v}}_h) \text{d}\bm{x}\right.\\
		&\left.\qquad\qquad + \sum_{e \in \partial_K \cap \Gamma}\int_{e}-(p^i+p_h^{s,N})\bm{n}\cdot(\overline{\bm{v}}-\overline{\bm{v}}_h) \text{d}s+\sum_{e \in \partial_K \cap \Gamma_s} \int_{e} \mathscr{T}^s_N \bm{u}_h^{N}(\overline{\bm{v}}-\overline{\bm{v}}_h)\text{d}s
		\right).
	\end{align*}
	Using the integration by parts, one can get
	\begin{align*}
		J_1&=\displaystyle \sum_{K \in \mathcal{M}_{h,a}}\left( \int_{K} 
		\left(\Delta p_h^{s,N}  +\kappa^2 p_h^{s,N}\right)(\overline{q}-\overline{q}_h) \text{d}\bm{x}
		-\displaystyle \sum_{e \in \partial_K}
		\int_{e} \nabla p_h^{s,N} \cdot \bm{\nu} (\overline{q}-\overline{q}_h)
		\text{d}s\right.\\
		&\left.\qquad\qquad + \sum_{e \in \partial_K \cap \Gamma}
		\int_{e}\left(  \partial_{\bm{n}} p^i 
		-\rho_f \omega^2 \bm{u}_h^N \cdot \bm{n}
		\right)(\overline{q}-\overline{q}_h)\text{d}s
		+  \sum_{e \in \partial_K \cap \Gamma_a} \int_{e} \mathscr{T}^a_N p_h^{s,N}(\overline{q}-\overline{q}_h)\text{d}s \right)\\
		&\quad+\displaystyle \sum_{K \in \mathcal{M}_{h,s}}\left(\int_{K}
		\left(\Delta^*\bm{u}_h^N+ \rho\omega^2\bm{u}_h^N \right)\cdot(\overline{\bm{v}}-\overline{\bm{v}}_h) \text{d}\bm{x}-\displaystyle \sum_{e \in \partial_K}
		\int_{e} \bm{Tu}_h^N \cdot (\overline{\bm{v}}-\overline{\bm{v}}_h) \text{d}s\right.\\
		&\left.\qquad\qquad\qquad + \displaystyle \sum_{e \in \partial_K \cap \Gamma}\int_{e}-(p^i+p_h^{s,N})\bm{n}\cdot(\overline{\bm{v}}-\overline{\bm{v}}_h) \text{d}s+\sum_{e \in \partial_K \cap \Gamma_s} \int_{e} \mathscr{T}^s_N \bm{u}_h^{N}(\overline{\bm{v}}-\overline{\bm{v}}_h)\text{d}s
		\right)\\
		&=\displaystyle \sum_{K \in \mathcal{M}_{h,a}}\left( \int_{K} 
		\mathcal{R}_a p_h^{s,N}(\overline{q}-\overline{q}_h) \text{d}\bm{x} +  \displaystyle \sum_{e \in \partial_K} \frac{1}{2}\int_{e}J_{e,a}(\overline{q}-\overline{q}_h)  \text{d}s \right)\\
		&\qquad+\displaystyle \sum_{K \in \mathcal{M}_{h,s}}\left(\int_{K}
		\mathcal{R}_s \bm{u}_h^{N} \cdot(\overline{\bm{v}}-\overline{\bm{v}}_h) \text{d}\bm{x}  +  \displaystyle \sum_{e \in \partial_K} \frac{1}{2}\int_{e}\bm{J_{e,s}}\cdot(\overline{\bm{v}}-\overline{\bm{v}}_h)  \text{d}s\right).
	\end{align*}	
	Let $q_h=\Pi_h^aq$ and $\bm{v}_h=\Pi_h^s\bm{v}$, where $\Pi_h^a$ and $\Pi_h^s$ are the Scott-Zhang \cite{Scott90}and Cl\'{e}ment \cite{Li22} interpolation operators, respectively, with the following interpolation estimates
	\begin{align*}
		&\|q-\Pi_h^aq \|_{L^2(K)} \lesssim h_{K} \|\nabla q\|_{L^2(\widetilde{K})},\quad
		\|q-\Pi_h^aq \|_{L^2(e)} \lesssim h_{e}^{1/2} \| q\|_{H^1(\widetilde{K}_e)},\\
		&\|\bm{v}-\Pi_h^s\bm{v} \|_{L^2(K)} \lesssim h_{K} \|\nabla \bm{v}\|_{L^2(\widetilde{K})},\quad
		\|\bm{v}-\Pi_h^s\bm{v} \|_{L^2(e)} \lesssim h_{e}^{1/2} \| \bm{v}\|_{H^1(\widetilde{K}_e)}.
	\end{align*}
	Here, $\widetilde{K}$ and $\widetilde{K}_e$
	are the union of all the elements in $\mathcal{M}_{h,s} \cup \mathcal{M}_{h,a}$, which have nonempty intersection with
	the element $K$ and the side $e$, respectively. 
	It follows from the Cauchy-Schwarz inequality and the above interpolation estimates that
	\begin{align*}
		|J_1| &\lesssim  \sum_{K \in \mathcal{M}_{h,a}} 
		\left(h_K\left\| \mathcal{R}_a p_h^{s,N}\right\|_{L^2(K)}
		\|\nabla q \|_{L^2(\widetilde{K})}+
		\sum_{e \in \partial_K} \frac{1}{2}h_e^{1/2}\| J_{e,a}\|_{L^2(e)} \|q\|_{H^1{(\widetilde{K}_e)}}	
		\right)\\
		&\quad +  \sum_{K \in \mathcal{M}_{h,s}}
		\left(h_K\left\| \mathcal{R}_s \bm{u}_h^{N}\right\|_{L^2(K)}
		\|\nabla\bm{v}\|_{L^2{(\widetilde{K})}}+
		\sum_{e \in \partial_K} \frac{1}{2}h_e^{1/2}\| \bm{J_{e,s}}\|_{L^2(e)}\|\bm{v}\|_{H^1{(\widetilde{K}_e)}}	
		\right)\\
		&\lesssim  \sum_{K \in \mathcal{M}_{h,a}} 
		\left(h_K\left\| \mathcal{R}_a p_h^{s,N}\right\|_{L^2(K)}+
		\left(  \sum_{e \in \partial_K}\frac{1}{2} h_e\| J_{e,a}\|_{L^2(e)}^2
		\right)^{1/2}
		\right)\|q\|_{H^1{(\Omega_a)}}\\
		&\quad +  \sum_{K \in \mathcal{M}_{h,s}}
		\left(h_K\left\| \mathcal{R}_s \bm{u}_h^{N}\right\|_{L^2(K)}+
		\left( \sum_{e \in \partial_K}  \frac{1}{2}h_e\| \bm{J_{e,s}}\|_{L^2(e)}^2
		\right)^{1/2}
		\right)\|\bm{v}\|_{H^1{(\Omega_s)}},
	\end{align*}
	Using (\ref{estimator1}) and (\ref{estimator2}), we have
	\begin{align}
		|J_1| &\lesssim \left(  \sum_{K \in \mathcal{M}_{h,a}} \eta_{K,a}^2 \right)^{1/2}\|q\|_{H^1{(\Omega_a)}}+
		\left(  \sum_{K \in \mathcal{M}_{h,s}} \eta_{K,s}^2 \right)^{1/2}\|\bm{v}\|_{H^1{(\Omega_s)}}\nonumber\\
		&\lesssim \left(  \sum_{K \in \mathcal{M}_{h,a}} \eta_{K,a}^2+
		\sum_{K \in \mathcal{M}_{h,s}} \eta_{K,s}^2 \right)^{1/2} \|\bm{V} \|_{\mathcal{H}^1(\Omega)}. \label{J_1}
	\end{align}
	
	The next work is to estimate $J_2$. It follows from (\ref{DtN1}) and (\ref{trDtN1}) that
	\begin{align*}
		|J_2|&=\left|\langle(\mathscr{T}^a-\mathscr{T}^a_N)p^s,q\rangle_{\Gamma_a}\right|=\left|\Lambda \sum_{|n| > N} \mathrm{i}\beta_np_n^{s}(b)\overline{q}_n(b)   \right|\\
		& \leq \Lambda \sum_{|n| > N} 
		\left|\beta_n\right| |p_n^{s}(b)| \left|q_n(b)\right|.
	\end{align*}
	By Lemma \ref{trace} and Lemma \ref{hatp}, we have
	\begin{align*}
		|J_2|& \lesssim \Lambda \sum_{|n| > N}e^{-\left(b-b^{\prime}\right)\left|\beta_n\right|} \left|\beta_n\right| \left|p^s_n(b^{\prime})\right|\left|q_n(b)\right|\\
		&\lesssim  \max_{|n| > N}\left(e^{-\left(b-b^{\prime}\right)\left|\beta_n\right|}\right)\left(\Lambda \sum_{\left|n\right|>N}\left|\beta_n\right|\left|p^s_n(b^{\prime})\right|^2\right)^{1 / 2}
		\left(\Lambda \sum_{\left|n\right|>N}\left|\beta_n\right|\left|q_n(b)\right|^2\right)^{1 / 2}\\
		&\lesssim \max_{|n| > N}\left(e^{-\left(b-b^{\prime}\right)\left|\beta_n\right|}\right) \|p^s\|_{H^{1 / 2}\left(\Gamma_a^{\prime}\right)} \|q\|_{H^{1 / 2}\left(\Gamma_a\right)}\\
		&\lesssim \max_{|n| > N}\left(e^{-\left(b-b^{\prime}\right)\left|\beta_n\right|}\right) \|p^s\|_{H^{1}\left(\Omega_a\right)} \|q\|_{H^{1}\left(\Omega_a\right)}.
	\end{align*}
	Using the stability estimate (\ref{regular1}), one can get 
	\begin{align}
		|J_2|& \lesssim \max_{|n| > N}\left(e^{-\left(b-b^{\prime}\right)\left|\beta_n\right|}\right) \|\bm{U} \|_{\mathcal{H}^1(\Omega)} \|\bm{V} \|_{\mathcal{H}^1(\Omega)}\nonumber\\
		&\lesssim \max_{|n| > N}\left(e^{-\left(b-b^{\prime}\right)\left|\beta_n\right|}\right)  \left(\left\|p^i \right\|_{L^2(\Gamma)}+\left\|\partial_{\bm{n}} p^i \right\|_{L^2(\Gamma)}\right)\|\bm{V} \|_{\mathcal{H}^1(\Omega)}. \label{J_2}
	\end{align}
	
	It remains to estimate $J_3$.
	A direct calculation yields that 
	\begin{align*}
		\phi_n^{(j)}(-b)=\phi_n^{(j)}\left(-b^{\prime}\right) e^{\mathrm{i} \beta_n^{(j)}\left(b-b^{\prime}\right)}.
	\end{align*} 
	By (\ref{algebraic equation}), we have
	\begin{align*}
		{\left[\begin{array}{l}
				u_n^{(1)}(-b) \\
				u_n^{(2)}(-b)
			\end{array}\right] } & =\frac{1}{\chi_n}\left[\begin{array}{cc}
			\mathrm{i} \alpha_n & -\mathrm{i} \beta_n^{(2)} \\
			-\mathrm{i} \beta_n^{(1)} & -\mathrm{i} \alpha_n
		\end{array}\right]\left[\begin{array}{cc}
			e^{\mathrm{i} \beta_n^{(1)}\left(b-b^{\prime}\right)} & 0 \\
			0 & e^{\mathrm{i} \beta_n^{(2)}\left(b-b^{\prime}\right)}
		\end{array}\right]\left[\begin{array}{cc}
			-\mathrm{i} \alpha_n & \mathrm{i} \beta_n^{(2)} \\
			\mathrm{i} \beta_n^{(1)} & \mathrm{i} \alpha_n
		\end{array}\right]\left[\begin{array}{l}
			u_n^{(1)}\left(-b^{\prime}\right) \\
			u_n^{(2)}\left(-b^{\prime}\right)
		\end{array}\right] \\
		& =P_n\left[\begin{array}{l}
			u_n^{(1)}\left(-b^{\prime}\right) \\
			u_n^{(2)}\left(-b^{\prime}\right)
		\end{array}\right],
	\end{align*}
	where
	\begin{align*}
		& P_n^{(11)}=\frac{1}{\chi_n}\left(\alpha_n^2 e^{\mathrm{i} \beta_n^{(1)}\left(b-b^{\prime}\right)}+\beta_n^{(1)} \beta_n^{(2)} e^{\mathrm{i} \beta_n^{(2)}\left(b-b^{\prime}\right)}\right), \\
		& P_n^{(12)}=\frac{\alpha_n \beta_n^{(2)}}{\chi_n}\left(-e^{\mathrm{i} \beta_n^{(1)}\left(b-b^{\prime}\right)}+e^{\mathrm{i} \beta_n^{(2)}\left(b-b^{\prime}\right)}\right), \\
		& P_n^{(21)}=\frac{\alpha_n \beta_n^{(1)}}{\chi_n}\left(-e^{\mathrm{i} \beta_n^{(1)}\left(b-b^{\prime}\right)}+e^{\mathrm{i} \beta_n^{(2)}\left(b-b^{\prime}\right)}\right), \\
		& P_n^{(22)}=\frac{1}{\chi_n}\left(\alpha_n^2 e^{\mathrm{i} \beta_n^{(2)}\left(b-b^{\prime}\right)}+\beta_n^{(1)} \beta_n^{(2)} e^{\mathrm{i} \beta_n^{(1)}\left(b-b^{\prime}\right)}\right) .
	\end{align*}
	It follows from (\ref{beta_j}) that $\beta_n^{(j)}$ is purely imaginary for sufficiently large $|n|$. By the mean value theorem, for sufficiently large $|n|$, there exists $\sigma \in\left(\mathrm{i} \beta_n^{(1)}, \mathrm{i} \beta_n^{(2)}\right)$ such that
	
	\begin{align*}
		\chi_n P_n^{(11)} & =\left(\alpha_n^2+\beta_n^{(1)} \beta_n^{(2)}\right) e^{\mathrm{i} \beta_n^{(1)}\left(b-b^{\prime}\right)}+\beta_n^{(1)} \beta_n^{(2)}\left(e^{\mathrm{i} \beta_n^{(2)}\left(b-b^{\prime}\right)}-e^{\mathrm{i} \beta_n^{(1)}\left(b-b^{\prime}\right)}\right) \\
		& =\left(\alpha_n^2+\beta_n^{(1)} \beta_n^{(2)}\right) e^{\mathrm{i} \beta_n^{(1)}\left(b-b^{\prime}\right)}+\beta_n^{(1)} \beta_n^{(2)}\left(b-b^{\prime}\right) \mathrm{i}\left(\beta_n^{(2)}-\beta_n^{(1)}\right) e^{\sigma\left(b-b^{\prime}\right)}.
	\end{align*}
	A simple calculation yields
	\begin{align*}
		\alpha_n^2+\beta_n^{(1)} \beta_n^{(2)} & =\alpha_n^2-\left(\alpha_n^2-\kappa_1^2\right)^{1 / 2}\left(\alpha_n^2-\kappa_2^2\right)^{1 / 2} \\
		& =\frac{\alpha_n^2\left(\kappa_1^2+\kappa_2^2\right)-\kappa_1^2 \kappa_2^2}{\alpha_n^2+\left(\alpha_n^2-\kappa_1^2\right)^{1 / 2}\left(\alpha_n^2-\kappa_2^2\right)^{1 / 2}}<\kappa_1^2+\kappa_2^2
	\end{align*}
	and
	\begin{align*}
		\mathrm{i} \beta_n^{(2)}-\mathrm{i} \beta_n^{(1)} & =\left(\alpha_n^2-\kappa_1^2\right)^{1 / 2}-\left(\alpha_n^2-\kappa_2^2\right)^{1 / 2} \\
		& =\frac{\kappa_2^2-\kappa_1^2}{\left(\alpha_n^2-\kappa_1^2\right)^{1 / 2}+\left(\alpha_n^2-\kappa_2^2\right)^{1 / 2}}<\frac{\kappa_2^2-\kappa_1^2}{2\left(\alpha_n^2-\kappa_2^2\right)^{1 / 2}},
	\end{align*}
	which give
	\begin{align*}
		\left|P_n^{(11)}\right| \lesssim e^{\mathrm{i} \beta_n^{(1)}\left(b-b^{\prime}\right)}+|n| e^{\sigma\left(b-b^{\prime}\right)} \lesssim|n| e^{\mathrm{i} \beta_n^{(2)}\left(b-b^{\prime}\right)} .
	\end{align*}
	Similarly, we have
	\begin{align*}
		\left|P_n^{(ij)}\right| \lesssim|n| e^{\mathrm{i} \beta_n^{(2)}\left(b-b^{\prime}\right)}, \quad i, j=1,2.
	\end{align*}
	Combining the above estimates leads to
	\begin{align*}
		\left|u_n^{(1)}(-b)\right|^2+\left|u_n^{(2)}(-b)\right|^2 \lesssim n^2 e^{2 \mathrm{i} \beta_n^{(2)}\left(b-b^{\prime}\right)}\left(\left|u_n^{(1)}\left(-b^{\prime}\right)\right|^2+\left|u_n^{(2)}\left(-b^{\prime}\right)\right|^2\right).
	\end{align*}
	By (\ref{DtN2}), (\ref{trDtN2}), Remark \ref{M^n bound} and Lemma \ref{trace}, we have
	\begin{align*}
		|J_3|&=\left|\langle(\mathscr{T}^s-\mathscr{T}^s_N)\bm{u},\bm{v}\rangle_{\Gamma_s}\right|\\
		&=\left|\Lambda \sum_{|n|>N}\left(M_n \boldsymbol{u}_n(-b)\right) \cdot \overline{\boldsymbol{v}}_n(-b)\right| \\
		& \lesssim \sum_{|n|>N}\left|\left(|n|^{\frac{1}{2}} \boldsymbol{u}_n(-b)\right) \cdot\left(|n|^{\frac{1}{2}} \overline{\boldsymbol{v}}_n(-b)\right)\right| \\
		& \lesssim\left(\sum_{|n|>N}|n|\left(\left|u_n^{(1)}(-b)\right|^2+\left|u_n^{(2)}(-b)\right|^2\right)\right)^{1 / 2}\left(\sum_{|n|>N}|n|\left(\left|v_n^{(1)}(-b)\right|^2+\left|v_n^{(2)}(-b)\right|^2\right)\right)^{1 / 2}\\
		& \lesssim\left(\sum_{|n|>N}|n|^3 e^{2 \mathrm{i} \beta_n^{(2)}\left(b-b^{\prime}\right)}\left(\left|u_n^{(1)}\left(-b^{\prime}\right)\right|^2+\left|u_n^{(2)}\left(-b^{\prime}\right)\right|^2\right)\right)^{1 / 2}\|\boldsymbol{v}\|_{H^{1 / 2}(\Gamma_s)} \\
		& \lesssim \max _{|n|>N}\left(|n| e^{\mathrm{i} \beta_n^{(2)}\left(b-b^{\prime}\right)}\right)\|\boldsymbol{u}\|_{H^{1 / 2}\left(\Gamma_s^{\prime}\right)}\|\boldsymbol{v}\|_{H^{1 / 2}(\Gamma_s)} \\
		& \lesssim \max _{|n|>N}\left(|n| e^{\mathrm{i} \beta_n^{(2)}\left(b-b^{\prime}\right)}\right)\|\boldsymbol{u}\|_{H^1(\Omega_s)}\|\boldsymbol{v}\|_{H^1(\Omega_s)} .
	\end{align*}
	Using the stability estimate (\ref{regular1}), one can get 
	\begin{align}
		\left|J_3\right| &\lesssim \max _{|n|>N}\left(|n| e^{\mathrm{i} \beta_n^{(2)}\left(b-b^{\prime}\right)}\right)\|\bm{U} \|_{\mathcal{H}^1(\Omega)} \|\bm{V} \|_{\mathcal{H}^1(\Omega)}\nonumber\\
		& \lesssim \max _{|n|>N}\left(|n| e^{- \left|\beta_n^{(2)}\right|\left(b-b^{\prime}\right)}\right)  \left(\left\|p^i \right\|_{L^2(\Gamma)}+\left\|\partial_{\bm{n}} p^i \right\|_{L^2(\Gamma)}\right)\|\bm{V} \|_{\mathcal{H}^1(\Omega)}.\label{J_3}
	\end{align}
	
	Collecting (\ref{J_1}), (\ref{J_2}) and (\ref{J_3}) yields
	\begin{align*}
		|J_1+J_2+J_3| \lesssim \left( \left( \displaystyle \sum_{K \in \mathcal{M}_{h,a}} \eta_{K,a}^2+
		\displaystyle \sum_{K \in \mathcal{M}_{h,s}} \eta_{K,s}^2 \right)^{1/2}+ \Theta \left(\left\|p^i \right\|_{L^2(\Gamma)}+\left\|\partial_{\bm{n}} p^i \right\|_{L^2(\Gamma)}\right)
		\right)\| \bm{V}\|_{\mathcal{H}^1(\Omega)},
	\end{align*}
	where 
	\begin{equation*}
		\Theta = \max\left(\max_{|n| > N}\left(e^{-\left(b-b^{\prime}\right)\left|\beta_n\right|}\right),\max _{|n|>N}\left(|n| e^{- \left(b-b^{\prime}\right)\left|\beta_n^{(2)}\right|}\right) \right).
	\end{equation*}
	Thus, we complete the proof.	
\end{proof}

Let $\hat{M}_n=-\frac{1}{2}\left(M_n+M_n^{\star}\right)$,  where $M_n^{\star}$ denotes the complex transpose of $M_n$. It follows from  (\ref{M^n}) that
\begin{align*}
	\hat{M}_n=-\frac{\mathrm{i}}{\chi_n}\left[\begin{array}{cc}
		\omega^2 \beta_n^{(1)} & -2\mu \alpha_n \chi_n+\omega^2 \alpha_n \\
		-\omega^2 \alpha_n+2\mu \alpha_n \chi_n &  \omega^2 \beta_n^{(2)}
	\end{array}\right].
\end{align*}
\begin{lemma}\label{M is positive}
	$\hat{M}_n$ is positive definite for sufficiently large $|n|$.
\end{lemma}
\begin{proof}
	It follows from (\ref{beta_j}) that $\beta_n^{(j)}$ is purely imaginary for sufficiently large $|n|$. 
	Since $\chi_n>\kappa_1^2>0$, one can get
	\begin{align*}
		\hat{M}_n^{(11)}=-\frac{\mathrm{i}}{\chi_n} \omega^2 \beta_n^{(1)}=\frac{\omega^2}{\chi_n}\left(\alpha_n^2-\kappa_1^2\right)^{1 / 2}>0.
	\end{align*}
	A simple calculation yields that
	\begin{align*}
		\chi_n^2 \operatorname{det} \hat{M}_n & =-\omega^4 \beta_n^{(1)} \beta_n^{(2)}-\left(2\mu \alpha_n \chi_n-\omega^2 \alpha_n\right)^2 \\
		& =-\mu^2 \kappa_2^4\left(\chi_n-\alpha_n^2\right)-\mu^2 \alpha_n^2\left(2\chi_n-\kappa_2^2\right)^2 \\
		& =\mu^2 \chi_n\left(-\kappa_2^4-4\alpha_n^2 \chi_n+4 \alpha_n^2 \kappa_2^2\right) .
	\end{align*}
	It follows from Lemma \ref{chi_n bound } that for sufficiently large $|n|$, there holds 
	\begin{align*}
		\operatorname{det}\hat{M}_n >0.
	\end{align*}
	Thus, the proof is completed.
\end{proof}

Based on Lemma \ref{M is positive}, one can prove the following estimate.
\begin{lemma}\label{Tn}
	Let $\Omega_s^{\prime}=\left\{\bm{x} \in \mathbb{R}^2:-b<x_2<-b^{\prime}, 0<x_1<\Lambda\right\}$. Then for any $\delta_1>0$, there exists a positive constant $C(\delta_1)$ independent of $N$, such that
	\begin{align*}
		\Re \int_{\Gamma_s} \mathscr{T}^s_N \boldsymbol{\zeta} \cdot \overline{\boldsymbol{\zeta}} \mathrm{d} s \leq C(\delta_1)\|\boldsymbol{\zeta}\|_{L^2\left(\Omega_s^{\prime}\right)}^2+\delta_1\|\boldsymbol{\zeta}\|_{H^1\left(\Omega_s^{\prime}\right)}^2.
	\end{align*}
\end{lemma}
\begin{proof}
	Using (\ref{trDtN2}), we have
	\begin{align*}
		\Re \int_{\Gamma_s} \mathscr{T}^s_N \boldsymbol{\zeta} \cdot \overline{\boldsymbol{\zeta}} \mathrm{d} s=\Lambda \sum_{|n| \leq N} \Re\left(M_n \boldsymbol{\zeta}_n\right) \cdot \overline{\boldsymbol{\zeta}_n}=-\Lambda \sum_{|n| \leq N}\left(\hat{M}_n \boldsymbol{\zeta}_n\right) \cdot \overline{\boldsymbol{\zeta}_n}.
	\end{align*}
	By Lemma \ref{M is positive}, $\hat{M}_n$ is positive definite for sufficiently large $|n|$. Hence, for fixed $\omega, \lambda, \mu$, there exists $N^*$ such that for any $n>N^*$, $-\left(\hat{M}_n \boldsymbol{\zeta}_n\right) \cdot \overline{\boldsymbol{\zeta}_n} \leq 0$ . Thus, $\Re \int_{\Gamma_s} \mathscr{T}^s_N \boldsymbol{\zeta} \cdot \bar{\boldsymbol{\zeta}}\text{d}s$ can be rewritten as
	\begin{align}\label{Tn zeta zeta}
		\Re \int_{\Gamma_s} \mathscr{T}^s_N \boldsymbol{\zeta} \cdot \overline{\boldsymbol{\zeta}} \mathrm{d} s=-\Lambda \sum_{|n| \leq \min \left(N^*, N\right)}\left(\hat{M}_n \boldsymbol{\zeta}_n\right) \cdot \overline{\boldsymbol{\zeta}_n}-\Lambda \sum_{N>|n|>\min \left(N^*, N\right)}\left(\hat{M}_n \boldsymbol{\zeta}_n\right) \cdot \overline{\boldsymbol{\zeta}_n}.
	\end{align}
	where $\sum_{N>|n|>\min \left(N^*, N\right)}\left(\hat{M}_n \boldsymbol{\zeta}_n\right) \cdot \overline{\boldsymbol{\zeta}_n}=0$ if $N>N^*$. Since the second term in the right hand side of (\ref{Tn zeta zeta}) is non-positive, we only need to estimate the first part in the right hand side of (\ref{Tn zeta zeta}). Hence, there exists a constant $C$ depending only on $\omega, \mu$ and $\lambda$, such that $\left|\left(\hat{M}_n \boldsymbol{\zeta}_n\right) \cdot \boldsymbol{\zeta}_n\right| \leq C\left|\boldsymbol{\zeta}_n\right|^2$ for all $|n| \leq \min \left(N^*, N\right)$.
	For any $\delta_1>0$, it follows from Yong's inequality that
	\begin{align*}
		\left(b-b^{\prime}\right)|\phi(-b)|^2&=\int_{-b}^{-b^{\prime}}|\phi(x_2)|^2 \mathrm{d} x_2+\int_{-b}^{-b^{\prime}} \int_{-b}^{x_2}\left(|\phi(s)|^2\right)^{\prime} \mathrm{d}s \mathrm{d} x_2 \\
		& \leq \int_{-b}^{-b^{\prime}}|\phi(x_2)|^2 \mathrm{d} x_2+\left(b-b^{\prime}\right) \int_{-b}^{-b^{\prime}} 2|\phi(x_2)|\left|\phi^{\prime}(x_2)\right| \mathrm{d} x_2 \\
		&=\int_{-b}^{-b^{\prime}}|\phi(  x_2)|^2 \mathrm{d} x_2+\left(b-b^{\prime}\right) \int_{-b}^{-b^{\prime}} 2 \frac{|\phi(  x_2)|}{\sqrt{\delta_1}} \sqrt{\delta_1}\left|\phi^{\prime}(x_2)\right| \mathrm{d} x_2 \\
		& \leq \int_{-b}^{-b^{\prime}}|\phi(  x_2)|^2 \mathrm{d} x_2+\frac{b-b^{\prime}}{\delta_1} \int_{-b}^{-b^{\prime}}|\phi(  x_2)|^2 \mathrm{d} x_2+\delta_1\left(b-b^{\prime}\right) \int_{-b}^{-b^{\prime}}\left|\phi^{\prime}(x_2)\right|^2 \mathrm{d} x_2,
	\end{align*}
	which gives
	\begin{align*}
		|\phi(-b)|^2 \leq\left[\frac{1}{\delta_1}+\left(b-b^{\prime}\right)^{-1}\right] \int_{-b}^{-b^{\prime}}|\phi(  x_2)|^2 \mathrm{d} x_2+\delta_1 \int_{-b}^{-b^{\prime}}\left|\phi^{\prime}(x_2)\right|^2 \mathrm{d} x_2 .
	\end{align*}
	Let $\phi(x_1,x_2)=\sum_{n \in \mathbb{Z}} \phi_n(x_2) e^{\mathrm{i} \alpha_n x_1}$. A straight calculation yields that
	\begin{align*}
		\|\nabla \phi\|_{L^2\left(\Omega_s^{\prime}\right)}^2 & =\Lambda \sum_{n \in \mathbb{Z}} \int_{-b}^{-b^{\prime}}\left(\left|\phi_n^{\prime}(x_2)\right|^2+\alpha_n^2\left|\phi_n(x_2)\right|^2\right) \mathrm{d} x_2, \\
		\|\phi\|_{L^2\left(\Omega_s^{\prime}\right)}^2 & =\Lambda \sum_{n \in \mathbb{Z}} \int_{-b}^{-b^{\prime}}\left|\phi_n(x_2)\right|^2 \mathrm{d} x_2 .
	\end{align*}
	It follows from the above estimates that for any $\phi \in H^1\left(\Omega_s^{\prime}\right)$, there holds
	\begin{align*}
		\|\phi\|_{L^2(\Gamma_s)}^2&=\Lambda \sum_{n \in \mathbb{Z}}\left|\phi_n(-b)\right|^2 \\
		& \leq \Lambda\left[\frac{1}{\delta_1}+\left(b-b^{\prime}\right)^{-1}\right] \sum_{n \in \mathbb{Z}} \int_{-b}^{-b^{\prime}}\left|\phi_n(x_2)\right|^2 \mathrm{d} x_2+\Lambda \delta_1 \sum_{n \in \mathbb{Z}} \int_{-b}^{-b^{\prime}}\left|\phi^{\prime}(x_2)\right|^2 \mathrm{d} x_2 \\
		& \leq \Lambda\left[\frac{1}{\delta_1}+\left(b-b^{\prime}\right)^{-1}\right] \sum_{n \in \mathbb{Z}} \int_{-b}^{-b^{\prime}}\left|\phi_n(x_2)\right|^2 \mathrm{d} x_2+\Lambda \delta_1 \sum_{n \in \mathbb{Z}} \int_{-b}^{-b^{\prime}}\left(\left|\phi_n^{\prime}(x_2)\right|^2+\alpha_n^2\left|\phi_n(x_2)\right|^2\right) \mathrm{d} x_2 \\
		& \leq\left[\frac{1}{\delta_1}+\left(b-b^{\prime}\right)^{-1}\right]\|\phi\|_{L^2\left(\Omega_s^{\prime}\right)}^2+\delta_1\|\nabla \phi\|_{L^2(\Omega_s^{\prime})}^2 \\
		& \leq C(\delta_1)\|\phi\|_{L^2\left(\Omega_s^{\prime}\right)}^2+\delta_1\|\nabla \phi\|_{L^2\left(\Omega_s^{\prime}\right)}^2 .
	\end{align*}
	Combining the above estimates, we are ready to get	
	\begin{align*}
		\Re \int_{\Gamma_s} \mathscr{T}^s_N \boldsymbol{\zeta} \cdot \overline{\boldsymbol{\zeta}} \mathrm{d} s & \leq C\|\boldsymbol{\zeta}\|_{L^2(\Gamma_s)}^2 \leq C(\delta_1)\|\boldsymbol{\zeta}\|_{L^2\left(\Omega_s^{\prime}\right)}^2+\delta_1 \int_{\Omega_s^{\prime}}|\nabla \boldsymbol{\zeta}|^2 \mathrm{d} \boldsymbol{x} \\
		& \leq C(\delta_1)\|\boldsymbol{\zeta}\|_{L^2\left(\Omega_s^{\prime}\right)}^2+\delta_1\|\boldsymbol{\zeta}\|_{H^1\left(\Omega_s^{\prime}\right)}^2,
	\end{align*}
	which completes the proof.
\end{proof}

Besides, we introduce the following lemma.
\begin{lemma}\label{errJandJN}
	Let $\bm{\Phi}=(\gamma,\bm{w})$ be the solution of the dual problem (\ref{dual}). We have
	\begin{align*}
		&\left|\langle(\mathscr{T}^a-\mathscr{T}^a_N)\xi, \gamma \rangle_{\Gamma_a} \right| \lesssim N^{-2} \|\xi \|_{H^1(\Omega_a)}^2,\\
		&\left|\langle(\mathscr{T}^s-\mathscr{T}^s_N)\bm{\zeta}, \bm{w} \rangle_{\Gamma_s} \right| \lesssim N^{-1} \|\bm{\zeta} \|_{H^1(\Omega_s)}^2.
	\end{align*}
\end{lemma}
\begin{proof}
	Since the proofs of the first and second estimates are essentially as the same as that for Lemma 4.5 in \cite{Wang15} and inequality  (5.18) in  \cite{Li20Yuan},
	we have omitted the details for brevity. Interested readers are referred to \cite{Li20Yuan,Wang15} for more details . 
\end{proof}

Now we are ready to give the main result.
\begin{theorem}\label{posterior estimate}
	Let $\bm{U}$ and $\bm{U_h^N}$ be the solutions of (\ref{varational}) and (\ref{varationalN}), respectively. Then, there exists a positive
	integer $N_0$ independent of $h$ such that for any
	$N > N_0$, the following a posteriori error estimate holds
	\begin{align*}
		\|\bm{U}-\bm{U_h^N}\|_{\mathcal{H}^1(\Omega)} \lesssim \left( \left( \displaystyle \sum_{K \in \mathcal{M}_{h,a}} \eta_{K,a}^2+
		\displaystyle \sum_{K \in \mathcal{M}_{h,s}} \eta_{K,s}^2 \right)^{1/2}+ \Theta\left(\left\|p^i \right\|_{L^2(\Gamma)}+\left\|\partial_{\bm{n}} p^i \right\|_{L^2(\Gamma)}\right)
		\right),
	\end{align*}
	where 
	\begin{equation*}
		\Theta = \max\left(\max_{|n| > N}\left(e^{-\left(b-b^{\prime}\right)\left|\beta_n\right|}\right),\max _{|n|>N}\left(|n| e^{- \left(b-b^{\prime}\right)\left|\beta_n^{(2)}\right|}\right) \right).
	\end{equation*}
\end{theorem}
\begin{proof}
	Following the subsection 4.4 in \cite{Wang15}, one can get 
	\begin{align}\label{Re}
		\Re \langle\mathscr{T}^a_N \xi,\xi\rangle_{\Gamma_a} \leq 0.
	\end{align}
	It follows from the Cauchy-Schwarz inequality and Lemma \ref{trace theorem} that for any $\delta_2>0$, there exists a positive constant $C(\delta_2)$ independent of $N$, such that
	\begin{align*}
		\left|\left(a_3(\bm{\zeta},\xi)+a_4(\xi,\bm{\zeta)} \right)\right| &\lesssim
		\|\bm{\zeta} \|_{L^2(\Gamma)}^2+\|\xi \|_{L^2(\Gamma)}^2\\
		&\lesssim
		\|\bm{\zeta} \|_{L^2(\Gamma)}^2+\|\xi \|_{L^2(\partial \Omega_a)}^2\\
		&\lesssim
		\|\bm{\zeta} \|_{L^2(\Omega_s)}\|\bm{\zeta} \|_{H^1(\Omega_s)}+\|\xi \|_{L^2( \Omega_a)}\|\xi \|_{H^1( \Omega_a)}\\
		&\lesssim
		C(\delta_2) \|\bm{\Psi} \|_{L^2(\Omega)}^2+\delta_2\|\bm{\Psi} \|_{\mathcal{H}^1(\Omega)}^2.
	\end{align*}
	By (\ref{H1}), (\ref{Re}), Lemma \ref{estimate} and Lemma \ref{Tn}, we have
	\begin{align*}
		\|\bm{\Psi}\|_{\mathcal{H}^1(\Omega)}^2 \leq& C_1 \left( \left( \displaystyle \sum_{K \in \mathcal{M}_{h,a}} \eta_{K,a}^2+
		\displaystyle \sum_{K \in \mathcal{M}_{h,s}} \eta_{K,s}^2 \right)^{1/2}+ \Theta \left(\left\|p^i \right\|_{L^2(\Gamma)}+\left\|\partial_{\bm{n}} p^i \right\|_{L^2(\Gamma)}\right)
		\right)\|\bm{\Psi}\|_{\mathcal{H}^1(\Omega)}\\
		&+\left(C(\delta_1)+C(\delta_2)+C_2 \right)\|\bm{\Psi}\|_{L^2(\Omega)}^2+C_3(\delta_1+\delta_2)\|\bm{\Psi} \|_{\mathcal{H}^1(\Omega)}^2,
	\end{align*}
	where $C_1>0$, $C_2>0$ and $C_3>0$ are constants independent of $h$ and $N$.
	It follows from (\ref{L2}), , Lemma \ref{errJandJN}, Lemma \ref{estimate} and (\ref{regular2}) that
	\begin{align*}
		\|\bm{\Psi}\|_{L^2(\Omega)}^2 \leq& C_4 \left( \left( \displaystyle \sum_{K \in \mathcal{M}_{h,a}} \eta_{K,a}^2+
		\displaystyle \sum_{K \in \mathcal{M}_{h,s}} \eta_{K,s}^2 \right)^{1/2}+ \Theta \left(\left\|p^i \right\|_{L^2(\Gamma)}+\left\|\partial_{\bm{n}} p^i \right\|_{L^2(\Gamma)}\right)
		\right)\|\bm{\Psi}\|_{L^2(\Omega)}\\
		&+\left(C_5N^{-2}+C_6N^{-1} \right)\|\bm{\Psi}\|_{\mathcal{H}^1(\Omega)}^2,
	\end{align*}
	where $C_4 > 0$, $C_5 > 0$ and $C_6 > 0$ are independent of $h$ and $N$. Combining the above two estimates,
	we have
	{\footnotesize
	\begin{align*}
		\|\bm{\Psi}\|_{\mathcal{H}^1(\Omega)}^2 &\leq \left(C_1+C_4 \left(C(\delta_1)+C(\delta_2)+C_2 \right)\right) \left( \left( \displaystyle \sum_{K \in \mathcal{M}_{h,a}} \eta_{K,a}^2+
		\displaystyle \sum_{K \in \mathcal{M}_{h,s}} \eta_{K,s}^2 \right)^{1/2}+ \Theta \left(\left\|p^i \right\|_{L^2(\Gamma)}+\left\|\partial_{\bm{n}} p^i \right\|_{L^2(\Gamma)}\right)
		\right)\|\bm{\Psi}\|_{\mathcal{H}^1(\Omega)}\\
		&+\left(\left(C(\delta_1)+C(\delta_2)+C_2 \right)\left(C_5N^{-2}+C_6N^{-1} \right)+C_3(\delta_1+\delta_2)\right)\|\bm{\Psi}\|_{\mathcal{H}^1(\Omega)}^2.
	\end{align*}}
	Choosing sufficiently large $N_0$ such that 
	\begin{align*}
		\left(\left(C(\delta_1)+C(\delta_2)+C_2 \right)\left(C_5N^{-2}+C_6N^{-1} \right)+C_3(\delta_1+\delta_2)\right) < 1/2,
	\end{align*}
	which completes the proof
	by taking $N > N_0$.
\end{proof}

\section{Implementation and Numerical Examples}\label{numerical}
In this section, we discuss the algorithmic implementation of the adaptive 
DtN-FEM and present numerical examples to demonstrate the effectiveness of
the proposed method.

\subsection{Adaptive Algorithm}

Based on the a posteriori error estimate from Theorem \ref{posterior estimate}, we use the PDE toolbox of MATLAB to
implement the adaptive $P_1$-FEM. It is shown
in Theorem \ref{posterior estimate} that the a posteriori error estimator consists of two parts: the finite
element discretization error $\epsilon_h$ and the DtN truncation error $\epsilon_N$, where
\begin{equation}\label{discretization error}
	\epsilon_h = \left( \displaystyle \sum_{K \in \mathcal{M}_{h,a}} \eta_{K,a}^2+
	\displaystyle \sum_{K \in \mathcal{M}_{h,s}} \eta_{K,s}^2 \right)^{1/2}, 
\end{equation}
\begin{equation}\label{truction error}
	\epsilon_N = \Theta \left(\left\|p^i \right\|_{L^2(\Gamma)}+\left\|\partial_{\bm{n}} p^i \right\|_{L^2(\Gamma)}\right).
\end{equation}
In the implementation, we choose the parameters $b^{\prime}$, $b$ and $N$ based on (\ref{truction error}) to ensure that  $\epsilon_N$ is small enough. For example, we can take $\epsilon_N \leq 10^{-8}$. For simplicity,
in the following numerical experiments, $b^{\prime}$ is chosen as $b^{\prime}=\max_{x_1 \in(0, \Lambda)} f(x_1)$, and $N$ is taken to be the smallest positive integer such that $\epsilon_N \leq 10^{-8}$
. The algorithm is shown in Algorithm \ref{algorithm} for the adaptive DtN-FEM
to solve problem
(\ref{varational}).  In our implementation, we always choose $\tau = 0.5$.

\begin{algorithm} 
	\caption{The DtN-AFEM algorithm for the acoustic-elastic interaction in periodic structures} 
	\label{algorithm}
	\begin{algorithmic}[1]
		\State Given the tolerance $\epsilon > 0, \tau \in (0,1)$ ; 
		\State Fix the computational domain $\Omega  $ by choosing $b$;
		\State Choose $b^{\prime}$ and $N$ such that $\epsilon_N \leq 10^{-8}$;
		\State Construct an initial triangulation $\mathcal{M}_h$ over $\Omega$ and compute error estimators;
		\While {$ \epsilon_h > \epsilon$}
		\State Refine the mesh $\mathcal{M}_h$ according to the strategy:
		
		\qquad if $\eta_{\hat{T}}> \tau \displaystyle \max_{T \in \mathcal{M}_h} \eta_T$, then refine the element $\hat{T} \in \mathcal{M}_h$;
		\State Denote the new mesh still by $\mathcal{M}_h$ and solve the discrete problem (\ref{varationalN}) on the new mesh $\mathcal{M}_h$;
		\State Compute the corresponding error estimators;
		\EndWhile 
	\end{algorithmic} 
\end{algorithm}

\subsection{Numerical Examples}

We present four numerical examples to demonstrate the effectiveness of the proposed method. The first example
concerns the scattering by a flat surface and has an exact solution; the second and third examples are constructed such that the
solutions have corner singularities; the last example has a curved fluid-solid interface.
\newline

\textbf{Example 1.} 
We first introduce a model problem for which analytical solutions are available
for the evaluation of accuracy; see \cite{Hu16}. Consider the simplest periodic structure, where the interface $\Gamma:=\left\{\left(x_1,x_2\right): x_2=0\right\}$ is a straight line. Choose the parameters  $\omega=1,\mu=1,\lambda=1,\rho=1, \rho_f=1$ and $\theta = \pi/6$. 
It can be verified that the exact solution can be written as
\begin{align*}
	p^{s}(\bm{x}) & =a_1 e^{\textup{i} \left(\alpha x_1+ \beta_0 x_2\right)}, \quad \bm{x} \in \Omega_a, \\
	\bm{u}(\bm{x}) & =a_2\left[\begin{array}{c}
		\alpha \\
		-\beta_0^{(1)}
	\end{array}\right] e^{\textup{i}\left(\alpha x_1-\beta_0^{(1)} x_2\right)}+a_3\left[\begin{array}{c}
		\beta_0^{(2)} \\
		\alpha
	\end{array}\right] e^{\textup{i}\left(\alpha x_1-\beta_0^{(2)} x_2\right)},\quad \bm{x} \in \Omega_s,
\end{align*}
where the coefficients $a_j \in \mathbb{C}, j=1,2,3$ are to be determined. It follows from the transmission conditions (\ref{transmission1}) and (\ref{transmission2}) that $a_j, j=1,2,3$ satisfy the following linear system
\begin{equation*}
	\left[\begin{array}{ccc}
		\textup{i} \beta_0 & \rho_f \omega^2 \beta_0^{(1)} & -\rho_f \omega^2 \alpha \\
		0 & 2 \textup{i} \mu \alpha \beta_0^{(1)} & 2 \textup{i} \mu  \left(\beta_0^{(2)}\right)^2-\textup{i} \mu \kappa_s^2 \\
		1 & 2 \textup{i} \mu \left(\beta_0^{(1)}\right)^2+\textup{i} \lambda \kappa_p^2 & -2 \textup{i} \mu \alpha \beta_0^{(2)}
	\end{array}\right]\left[\begin{array}{l}
		a_1 \\
		a_2 \\
		a_3
	\end{array}\right]=\left[\begin{array}{c}
		\textup{i} \beta_0 \\
		0 \\
		-1
	\end{array}\right].
\end{equation*}
Take the period $\Lambda=4$.
\begin{figure}
	\centering
	\subfigure[$\Re p^s$]{\includegraphics[width=0.32\textwidth]{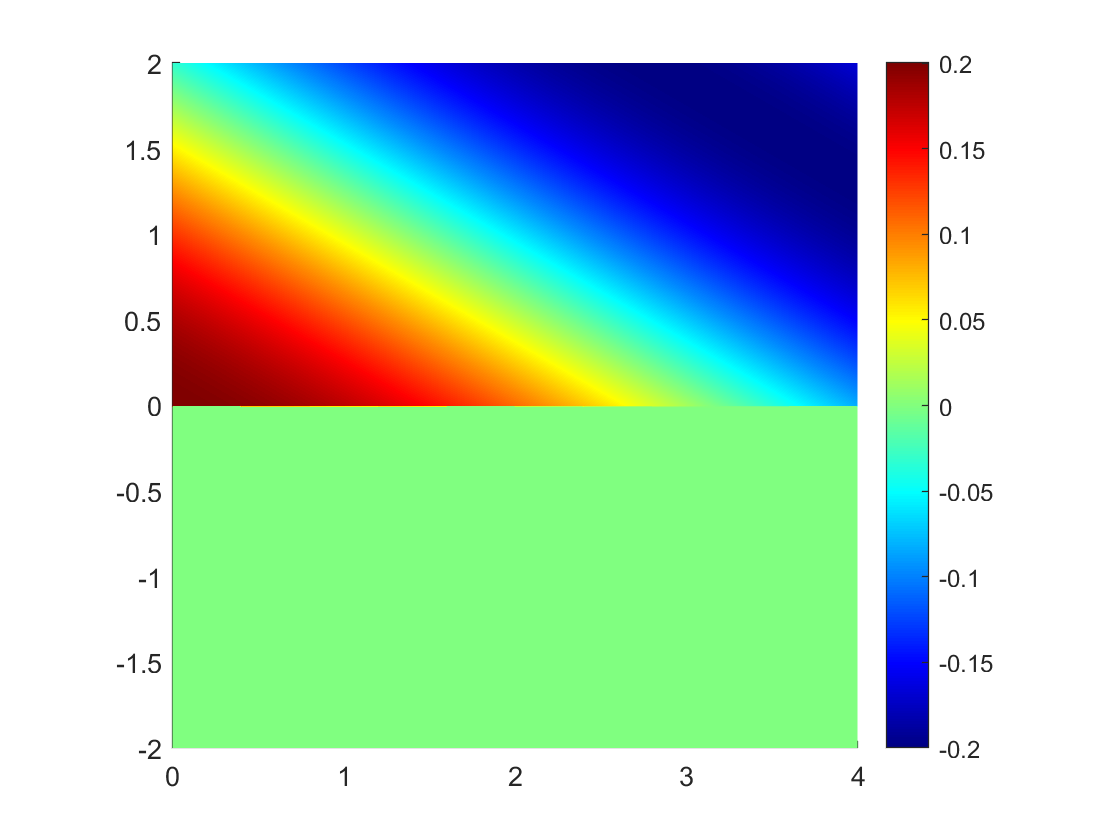}} 
	\subfigure[$\Re u_1$]{\includegraphics[width=0.32\textwidth]{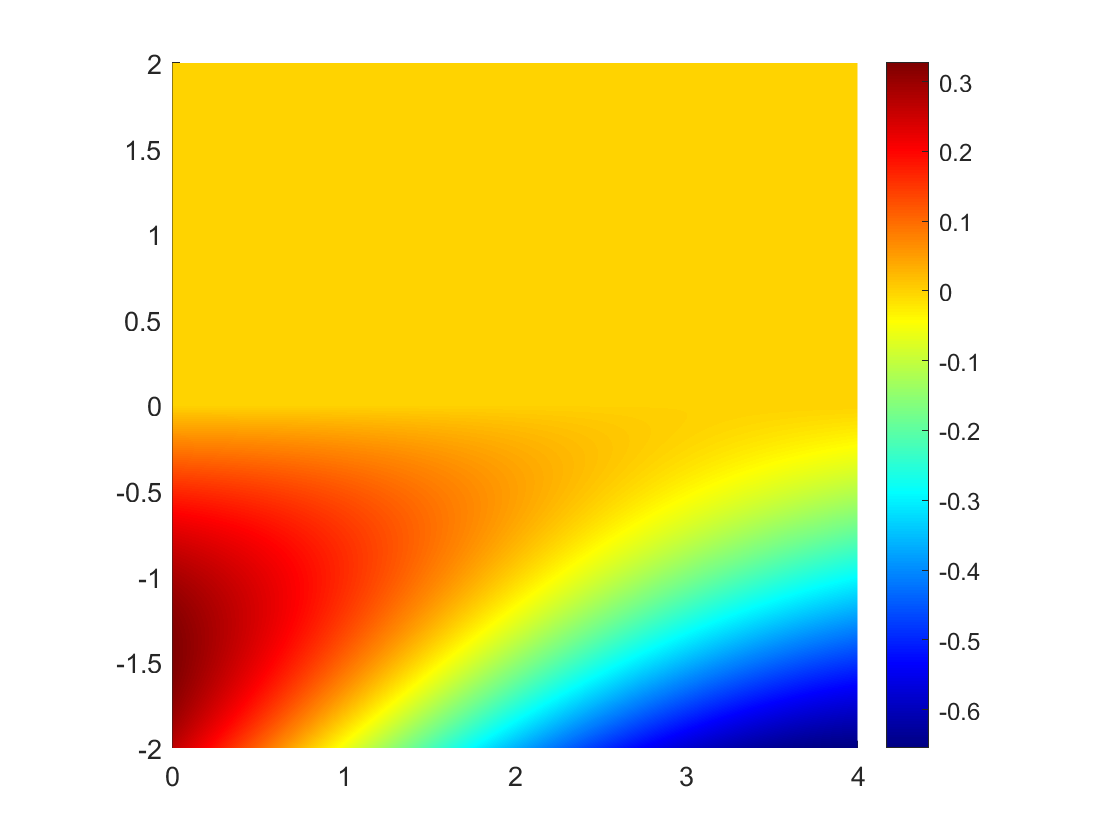}} 
	\subfigure[$\Re u_2$]{\includegraphics[width=0.32\textwidth]{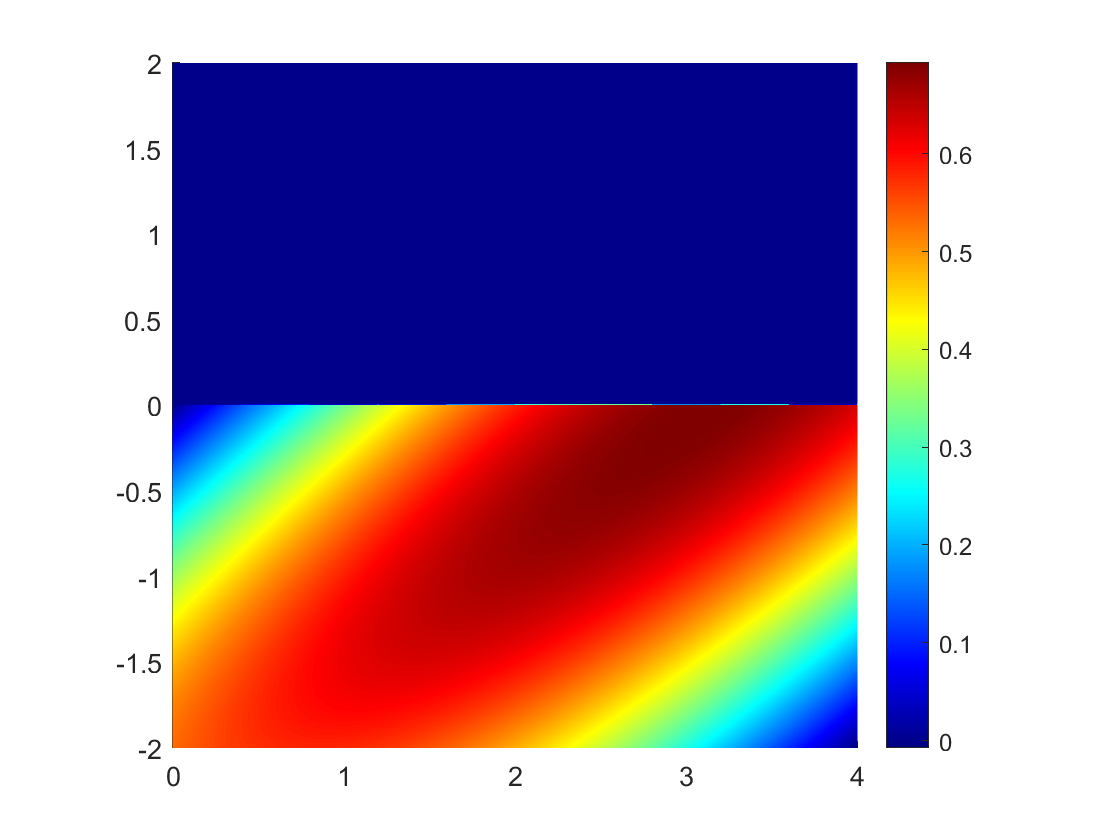}} \\
	\subfigure[$\Re p_h^{s,N}$]{\includegraphics[width=0.32\textwidth]{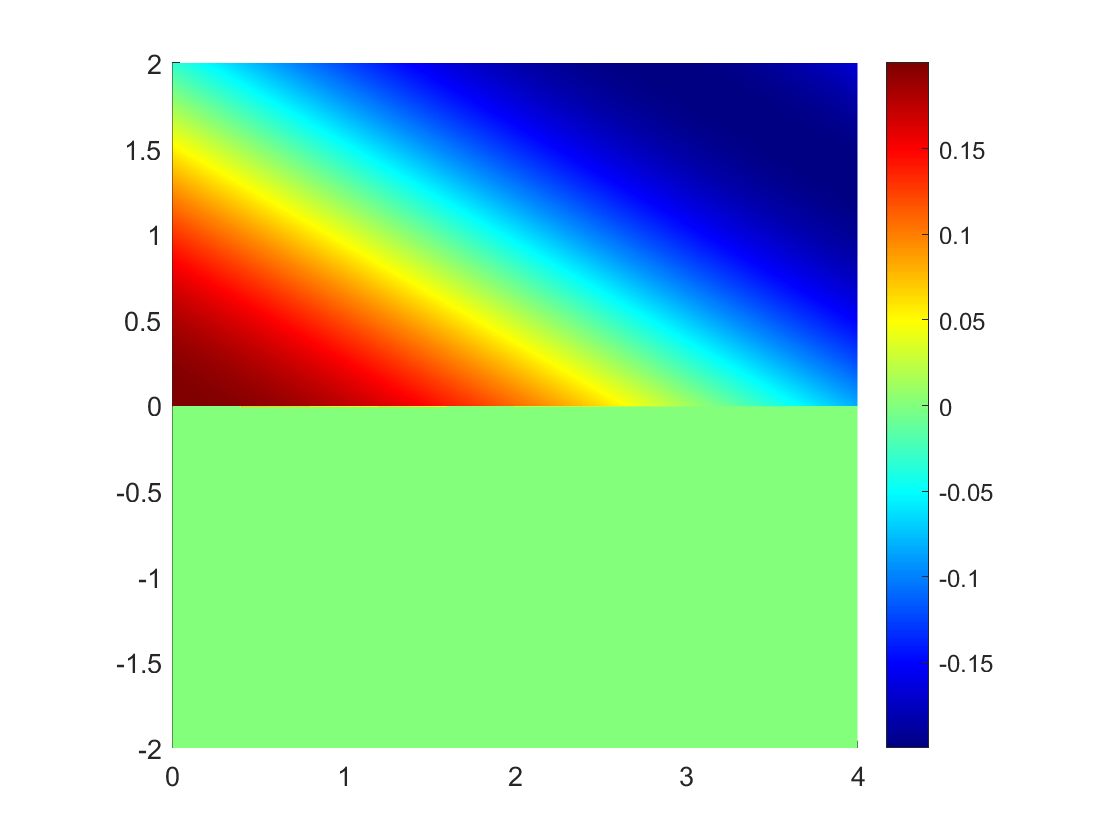}} 
	\subfigure[$\Re u_{h,1}^{N}$]{\includegraphics[width=0.32\textwidth]{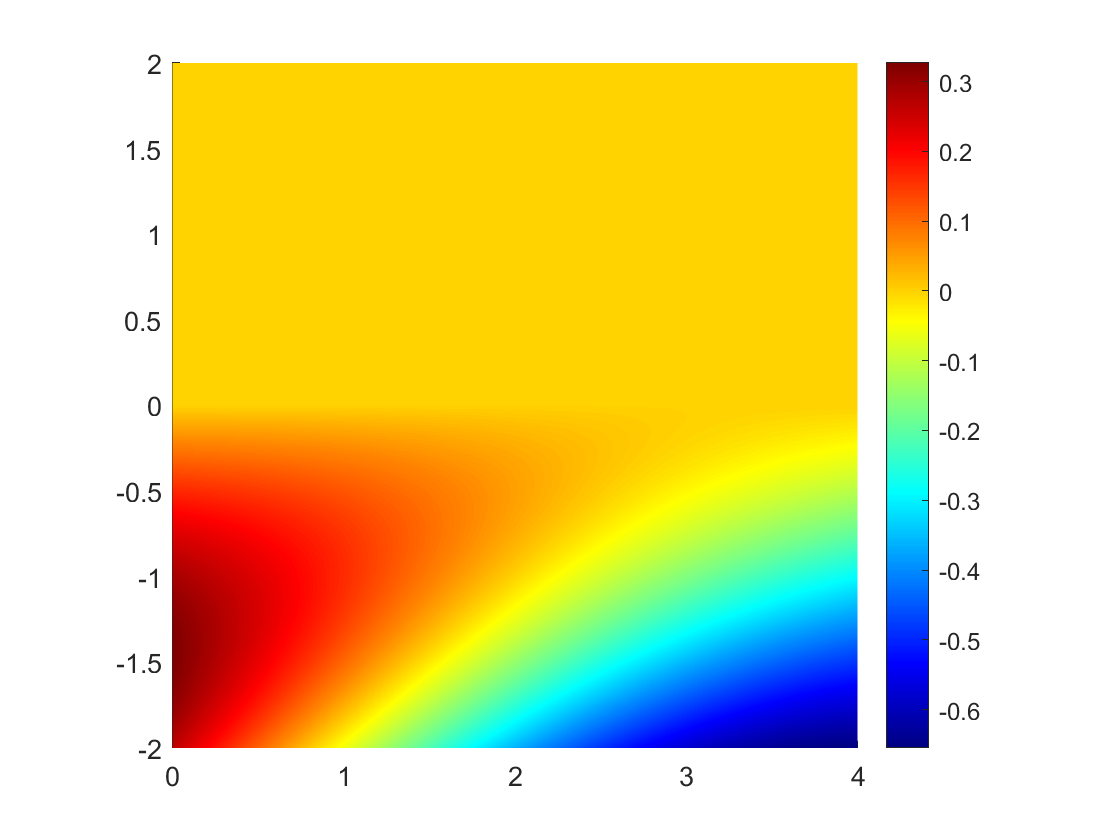}} 
	\subfigure[$\Re u_{h,2}^{N}$]{\includegraphics[width=0.32\textwidth]{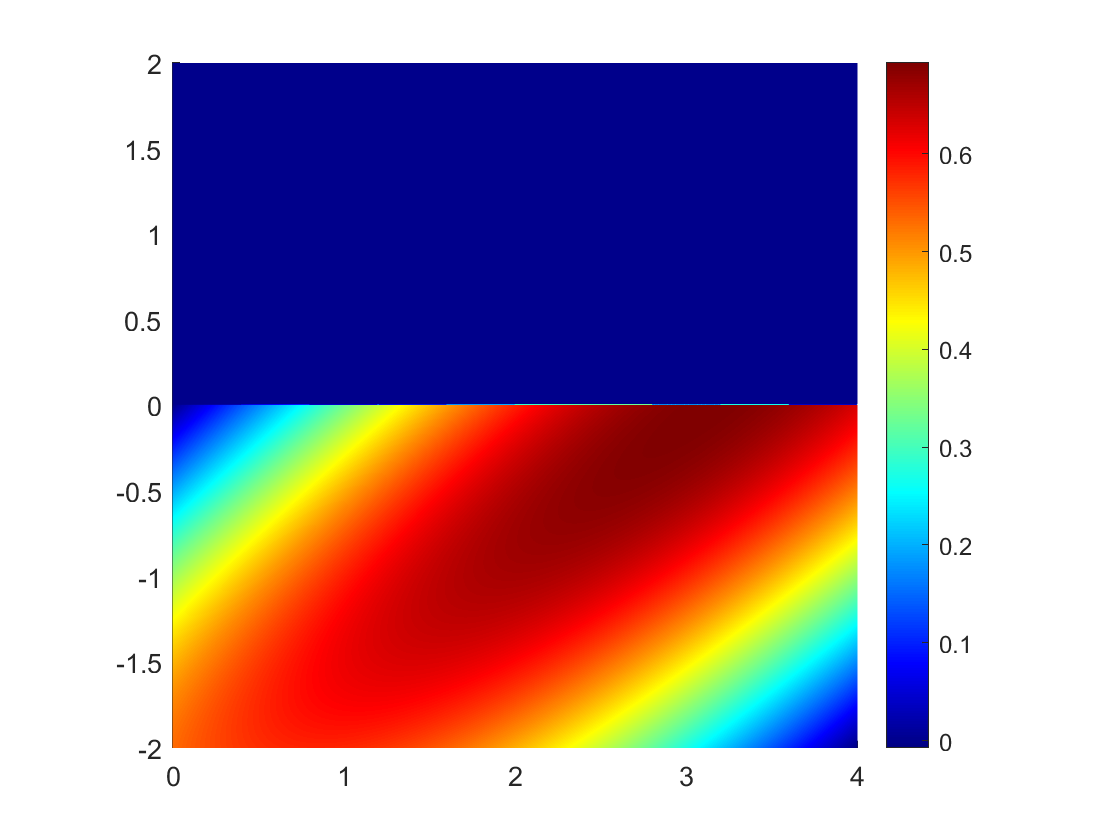}}\\
	\caption{The real parts of the exact solutions (a,b,c) and numerical solutions (d,e,f) with $\kappa = 1$ in Example 1}
	\label{contour of solution k=1}
	\vspace{0.2in}
\end{figure}

\begin{figure}
	\centering 
	\subfigure[A priori error ]{\includegraphics[width=0.49\textwidth]{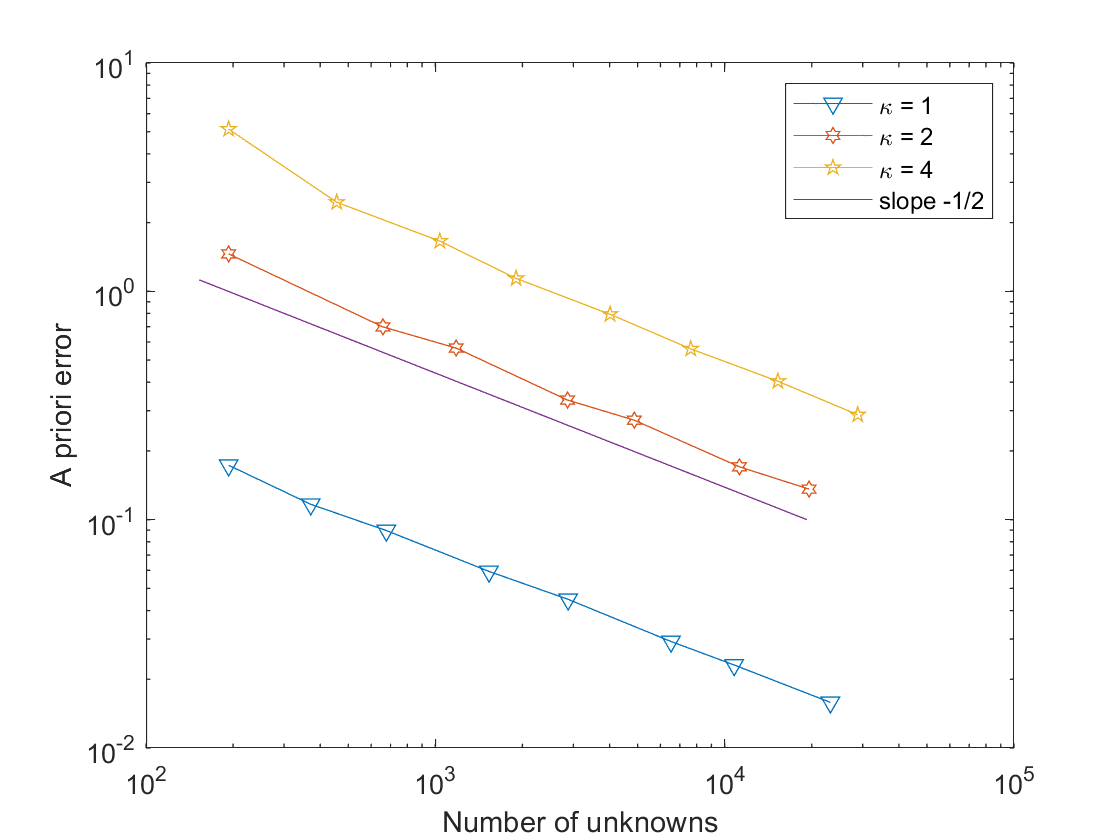}} 
	\subfigure[A posterior error]{\includegraphics[width=0.49\textwidth]{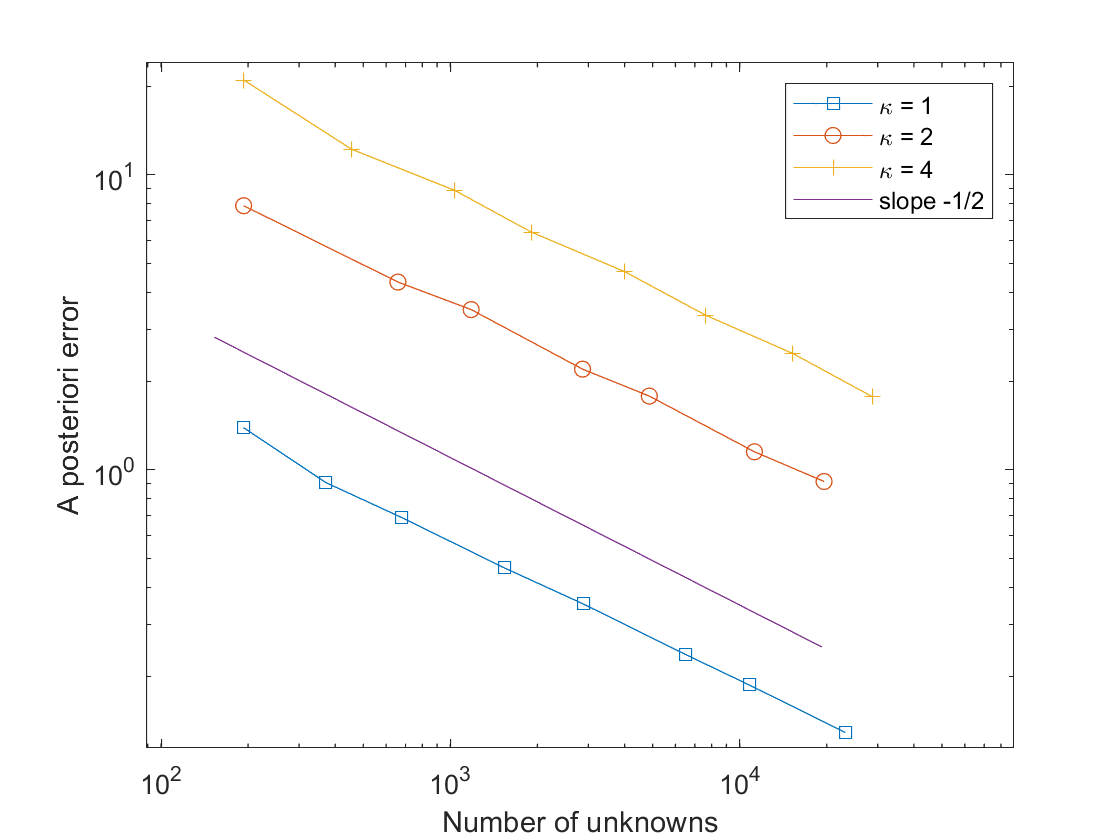}} 
	\caption{Quasi-optimality of the a priori and a posteriori errors in Example 1}
	\label{mesh}
\end{figure}

Denote the a priori error by $e_h=\left\|\boldsymbol{U}-\boldsymbol{U}_h^N\right\|_{\mathcal{H}^1(\Omega)}$.
We firstly check the accuracy of the adaptive DtN-FEM.
Fig. \ref{contour of solution k=1} shows the real parts of the exact and numerical solutions with $\kappa = 1$. It can be
observed that the numerical solutions are in a perfect agreement with the exact solutions. Fig. \ref{mesh} shows
the curves of $\text{log} \epsilon_h$ and $\text{log} e_h$ against $\text{log DoF}_h$ with $\kappa=1,~ 2,  ~4$, respectively,
where $\text{DoF}_h$
denotes the number of nodal points of the mesh $\mathcal{M}_h$. It proves that the decays of the a priori and a posteriori errors are both quasi-optimal.
\newline

\begin{figure}
	\centering
	\subfigure[$\Re p_h^{s,N}$]{\includegraphics[width=0.32\textwidth]{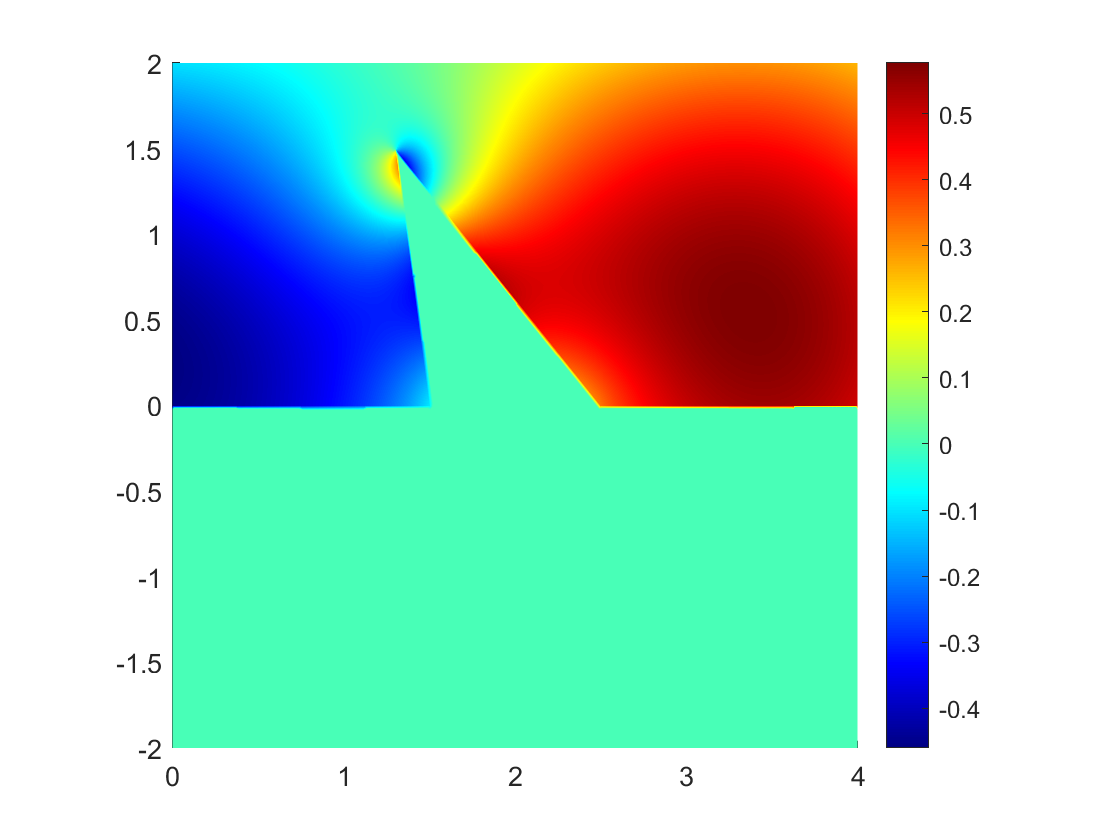}}
	\subfigure[$\Re u_{h,1}^{N}$]{\includegraphics[width=0.32\textwidth]{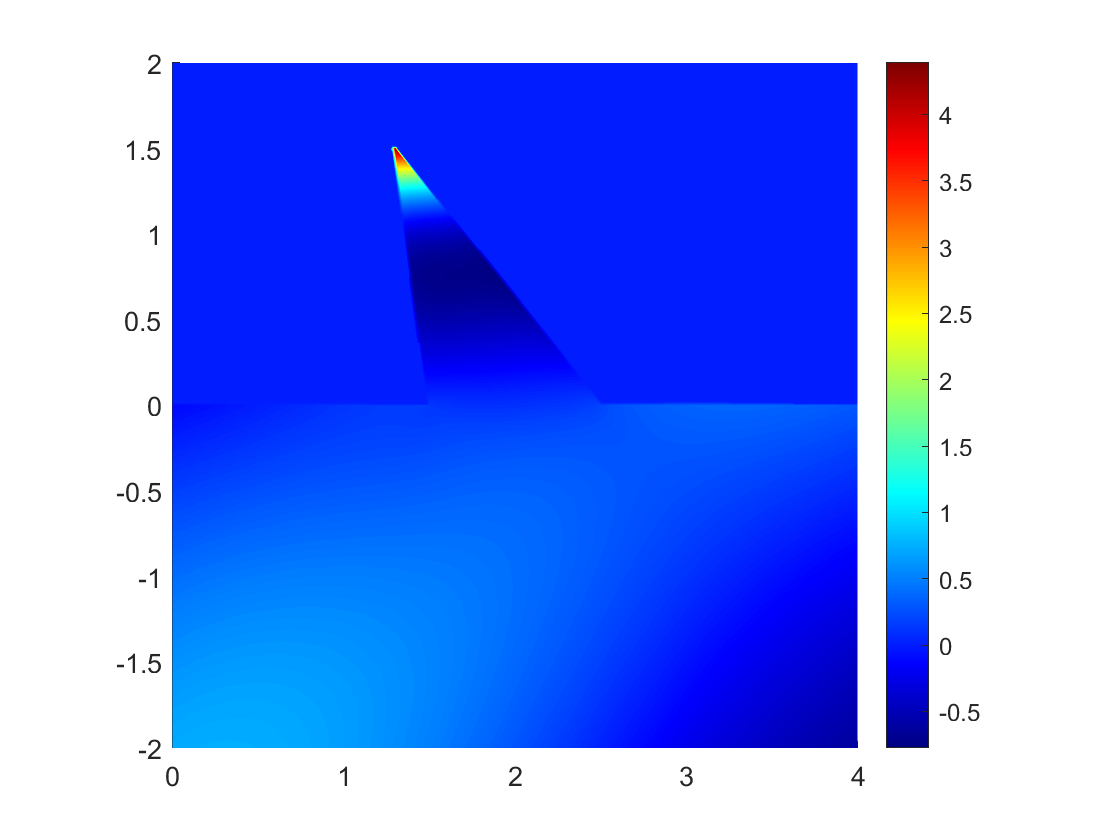}} 
	\subfigure[$\Re u_{h,2}^{N}$]{\includegraphics[width=0.32\textwidth]{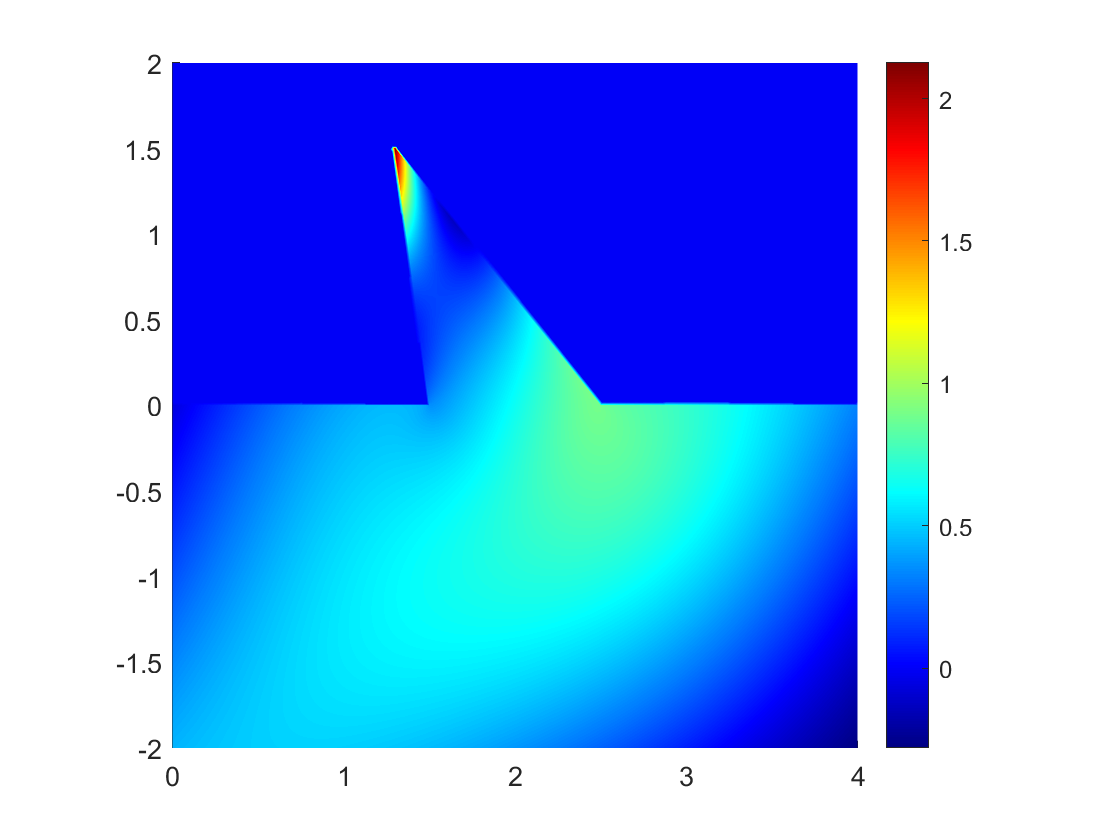}} 
	\caption{The real parts of the numerical solutions with $\kappa =1$ in Example 2}
	\label{k=1 example 2}
\end{figure}

\begin{figure}
	\centering
	\subfigure[A posterior error]{\includegraphics[width=0.495\textwidth]{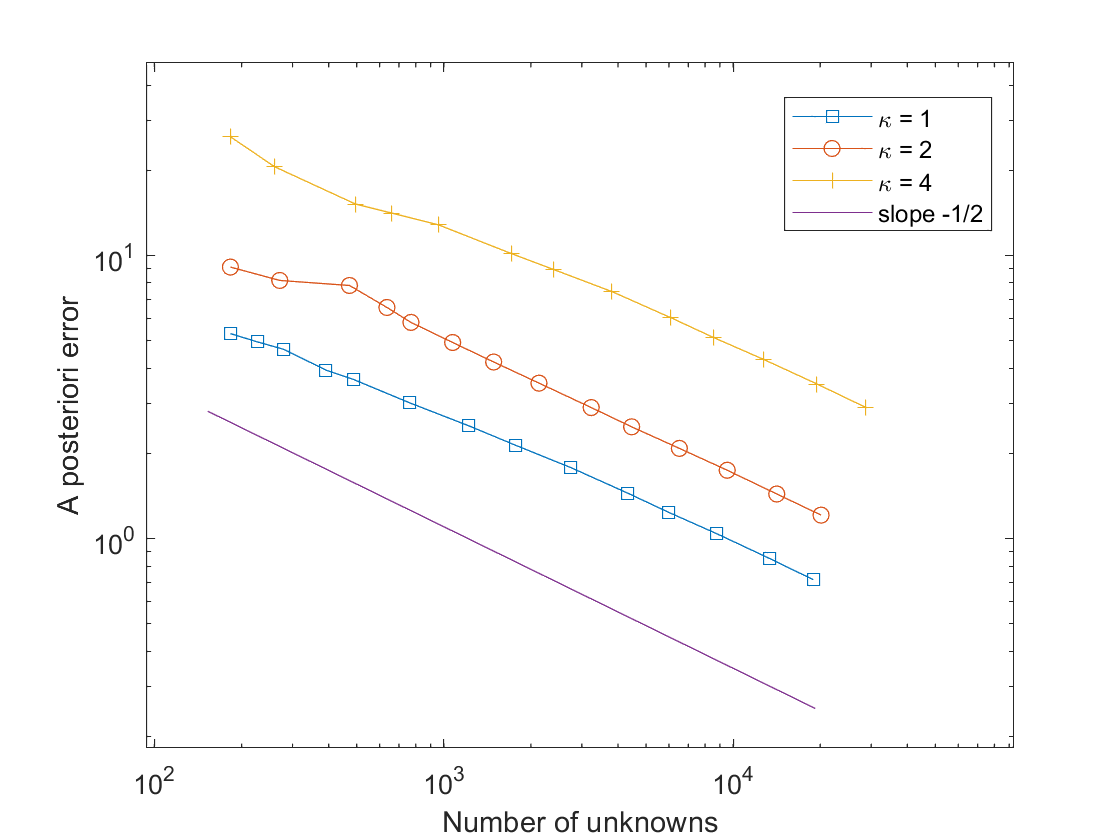}}
	\subfigure[Adaptive mesh]{\includegraphics[width=0.495\textwidth]{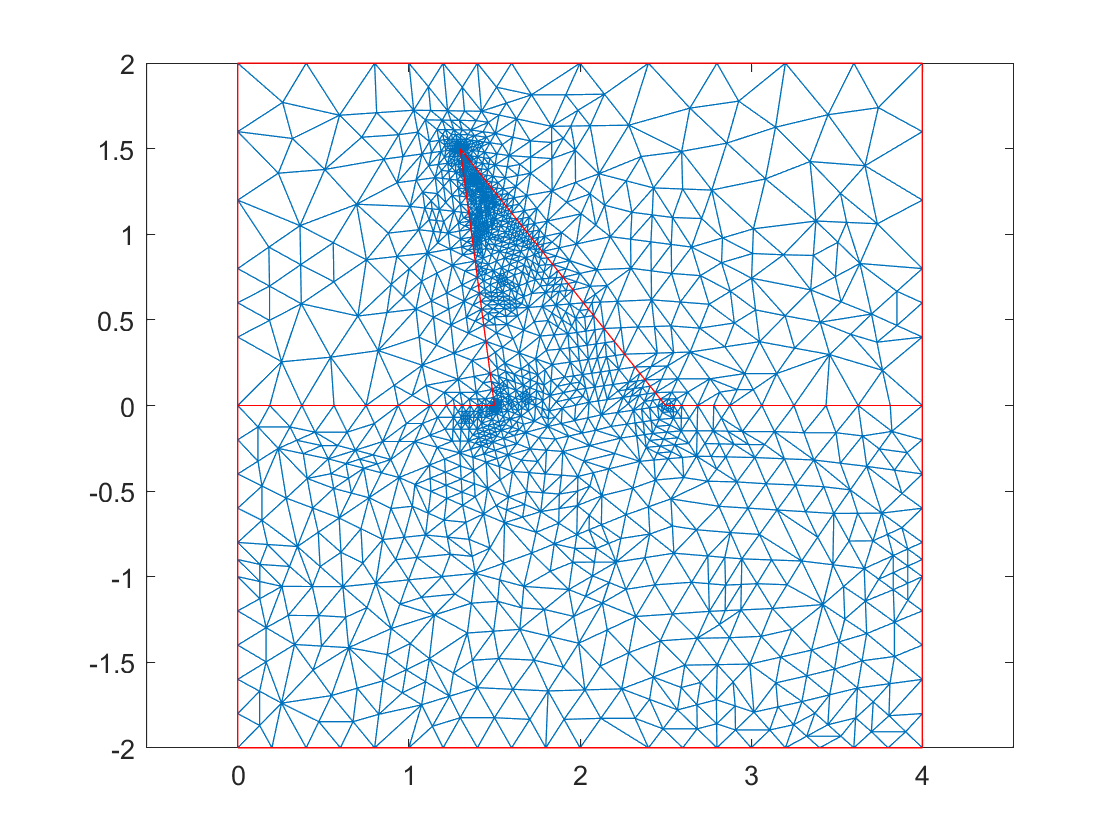}}
	\caption{Quasi-optimality of the a posteriori errors with different $\kappa$ (a); the adaptive mesh with 5,966 elements after 8 refinement iterations (b) in Example 2}
	\label{posterior_error and mesh example2}
\end{figure}

\begin{figure}
	\centering
	\subfigure[$\Re p_h^{s,N}$]{\includegraphics[width=0.32\textwidth]{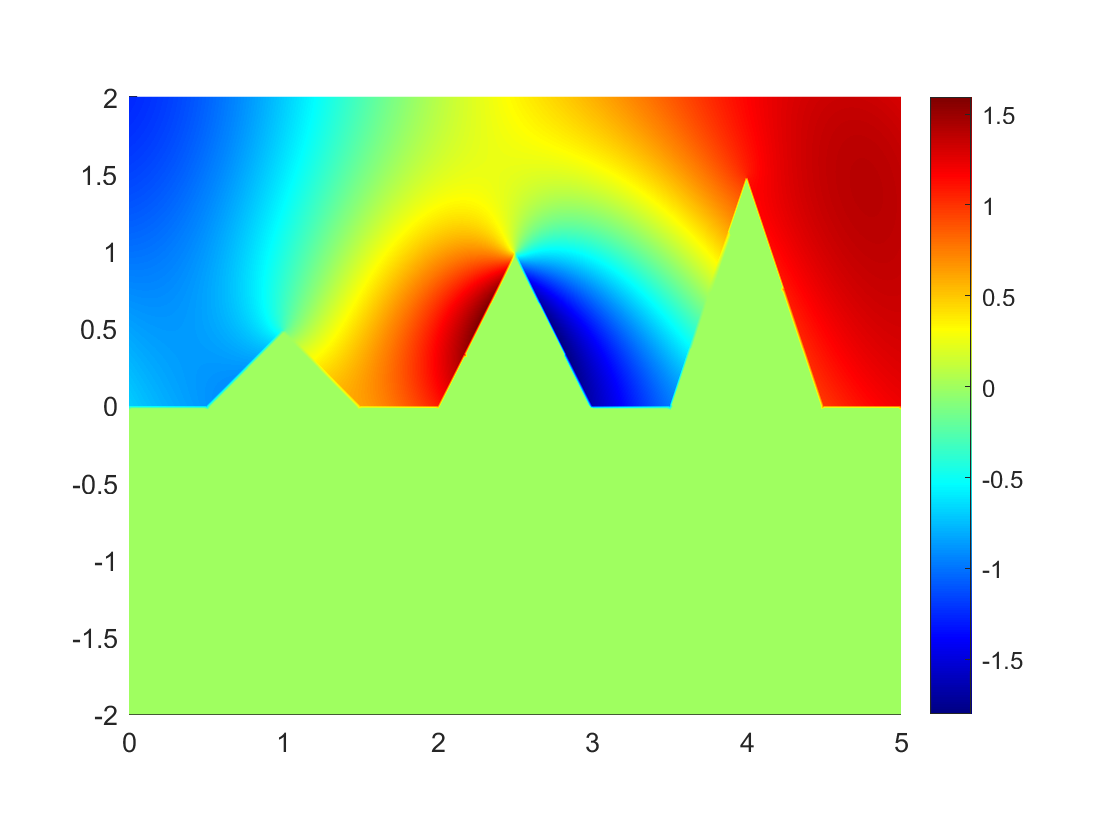}}
	\subfigure[$\Re u_{h,1}^{N}$]{\includegraphics[width=0.32\textwidth]{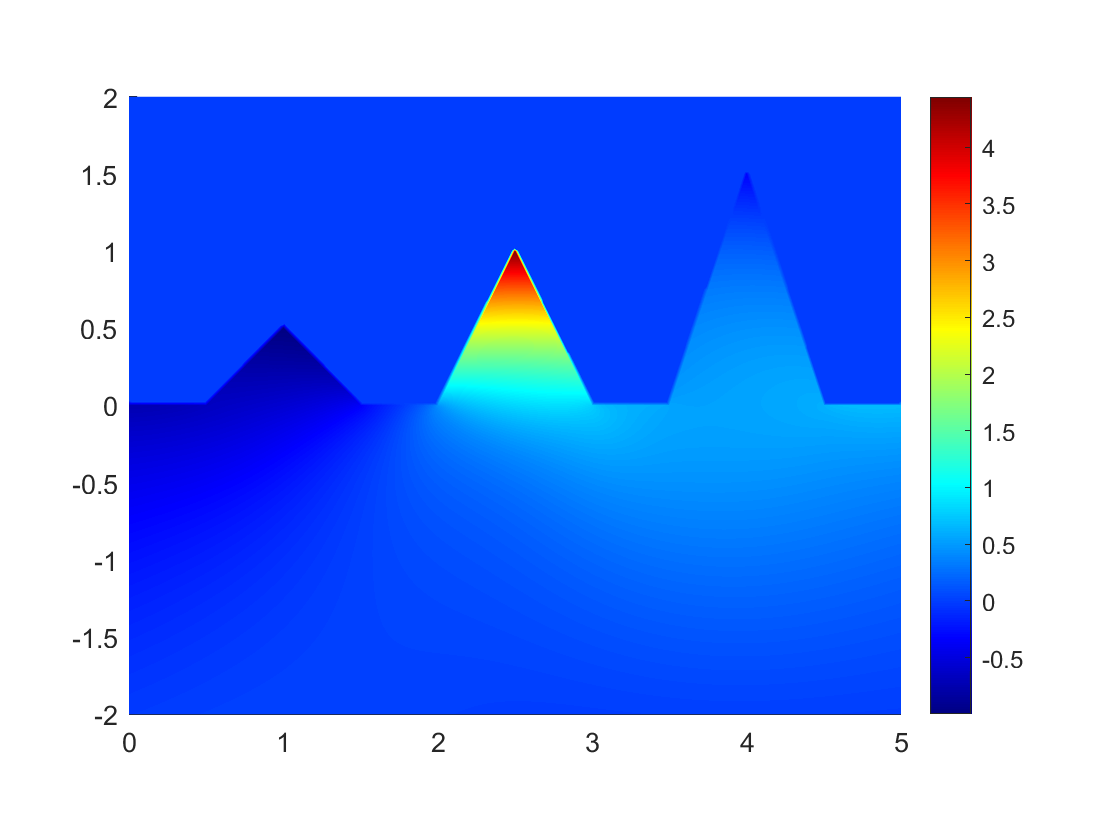}} 
	\subfigure[$\Re u_{h,2}^{N}$]{\includegraphics[width=0.32\textwidth]{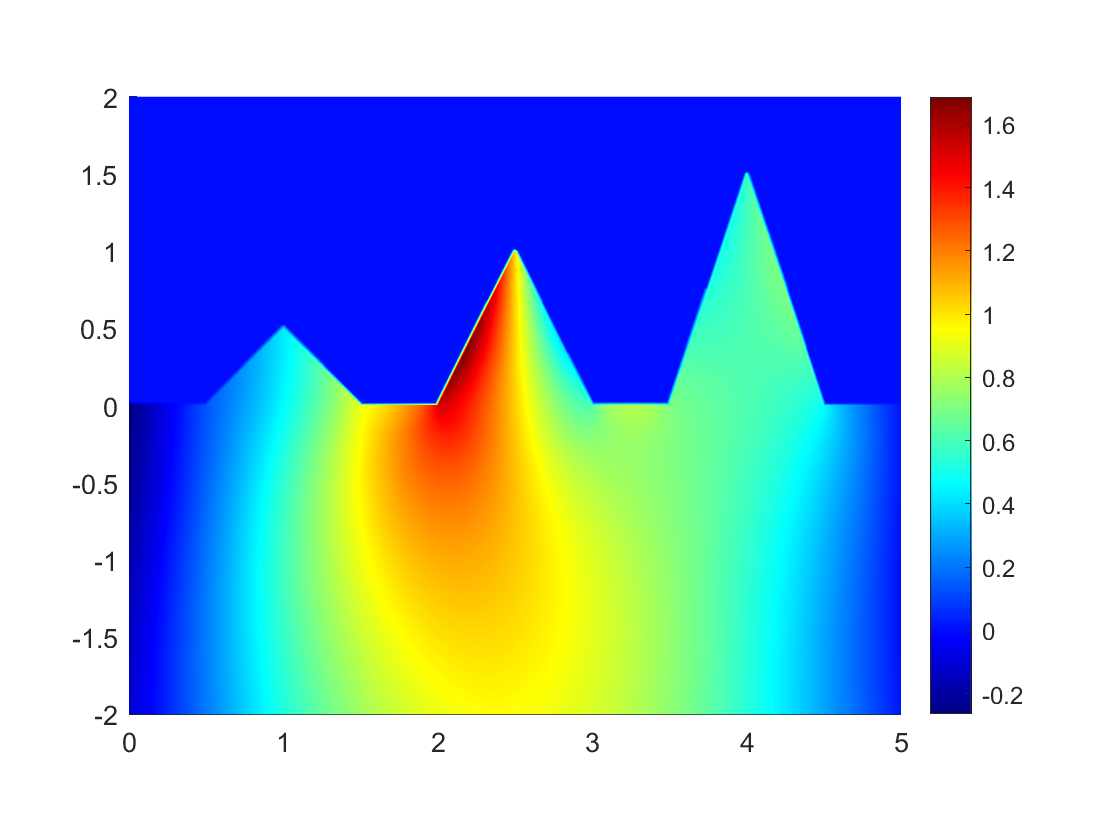}} 
	\caption{The real parts of the numerical solutions with $\kappa =1$ in Example 3}
	\label{k=1 example 3}
\end{figure}

\begin{figure}
	\centering
	\subfigure[A posterior error]{\includegraphics[width=0.495\textwidth]{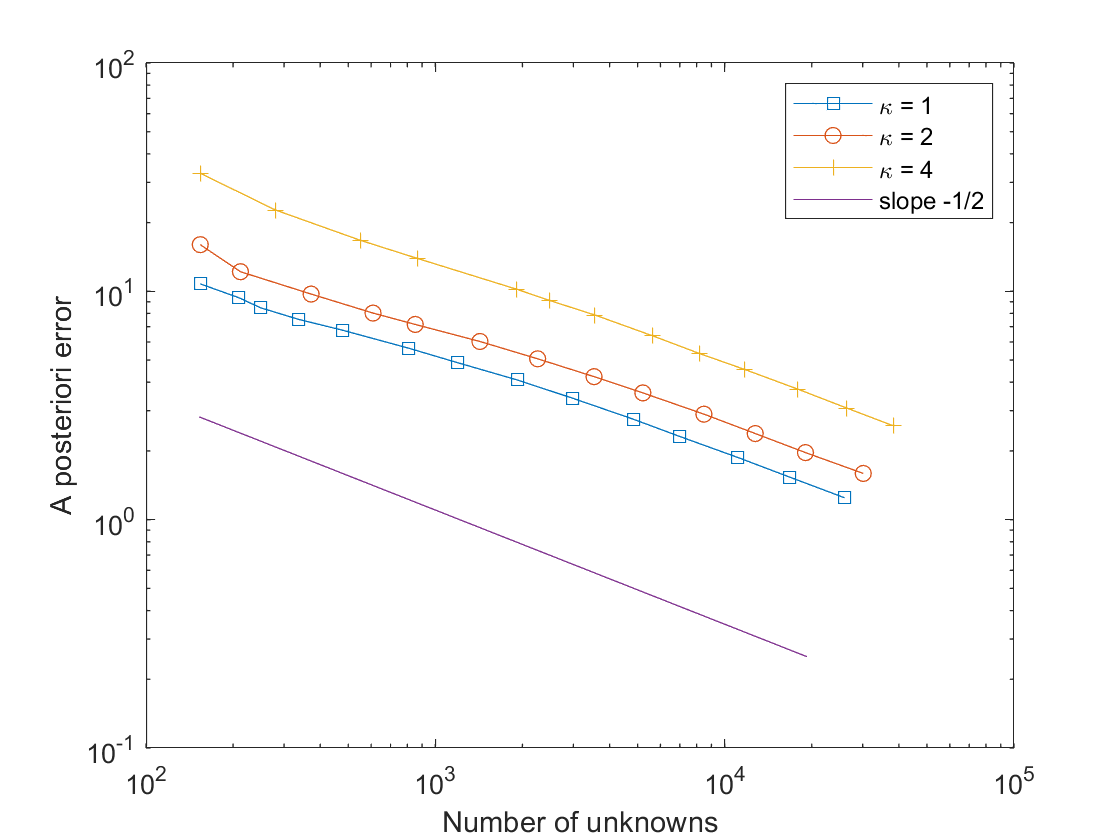}}
	\subfigure[Adaptive mesh]{\includegraphics[width=0.495\textwidth]{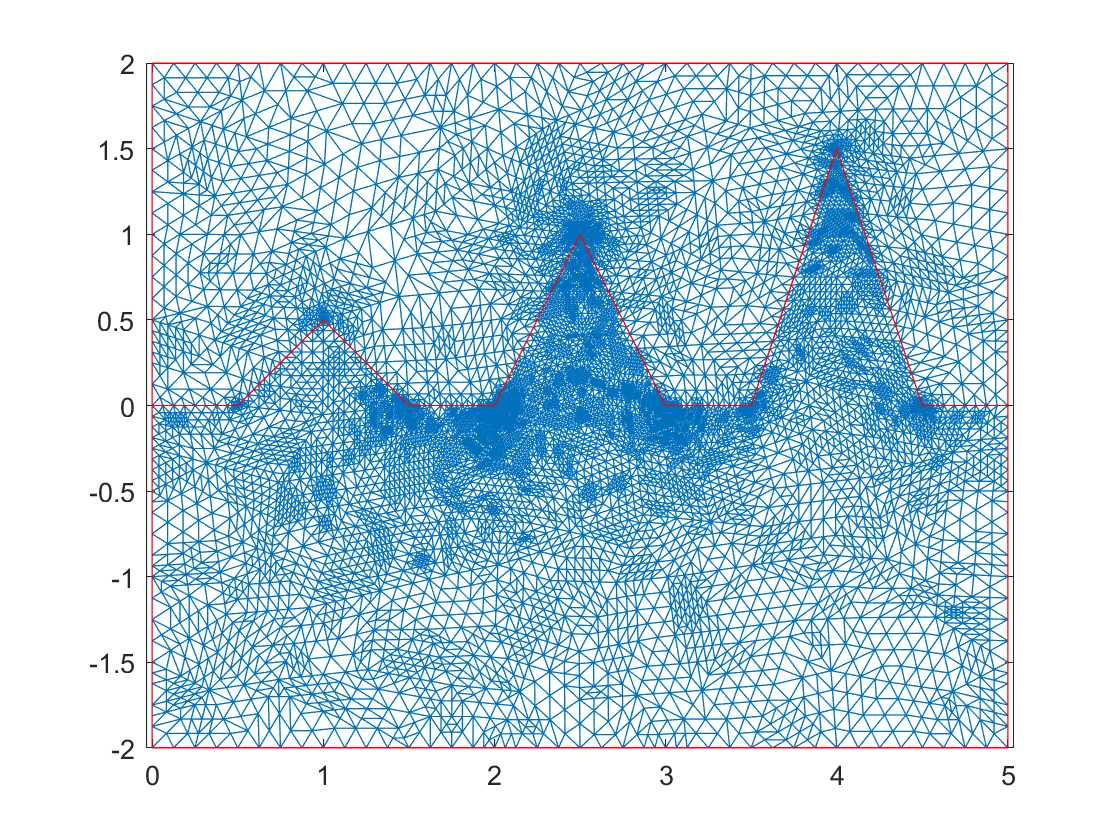}}
	\caption{Quasi-optimality of the a posteriori errors with different $\kappa$ (a); the adaptive mesh with 11,168 elements after 11 refinement iterations (b) in Example 3}
	\label{posterior_error and mesh}
\end{figure}

\begin{table}
	\centering
	\caption{Comparison of numerical results using uniform and adaptive refinements with  $\kappa = 1$ in Example 2}
	\begin{tabular}{lllll}
		\hline \multicolumn{2}{l}{ Uniform mesh } && \multicolumn{2}{l}{Adaptive  mesh } \\
		\hline
		 $\mathrm{DoF}_h$ & $\epsilon_h$  && $\mathrm{DoF}_h$  & $\epsilon_h$   \\
		\hline
		183 & 5.2990 && 183 & 5.2990  \\
		689 & 4.4728 && 391 & 3.9391 \\
		2673 & 3.5372 && 764 & 3.0137  \\
		10529 & 2.4110 && 1774 & 2.1346  \\
		41793 & 1.5043 && 4317 & 1.4482  \\
		\hline
	\end{tabular}
	\label{example2table}
\end{table}

\begin{table}
	\centering
	\caption{Comparison of numerical results using uniform and adaptive refinements with  $\kappa = 1$ in Example 3}
	\begin{tabular}{lllll}
		\hline \multicolumn{2}{l}{ Uniform mesh } && \multicolumn{2}{l}{Adaptive  mesh } \\
		\hline
		 $\mathrm{DoF}_h$ & $\epsilon_h$  && $\mathrm{DoF}_h$  & $\epsilon_h$   \\
		\hline
		154 & 10.7777 && 154 & 10.7777  \\
		577 & 8.3513 && 326 & 7.5627 \\
		2233 & 5.6889 && 808 & 5.6158  \\
		8785 & 3.6623 && 2984 & 3.3983  \\
		34849 & 2.3395 && 6953 & 2.3224  \\
		\hline
	\end{tabular}
	\label{example3table}
\end{table}

\textbf{Example 2.} This example concerns the scattering of the plane wave 
by a grating surface
with a sharp angle. The incident angle is chosen as $\theta = \pi/4$ and the remain parameters are the same as those
in Example 1. This example does not have an analytical solution.

The real parts of the numerical solutions with $\kappa = 1$ are presented in Fig \ref{k=1 example 2}. Fig. \ref{posterior_error and mesh example2} (a) presents the curves of $\log \epsilon_h$ versus $\log \mathrm{DoF}_h$ with different wave number $\kappa$, where the decays of the a posteriori errors are $\mathcal{O}\left(\mathrm{DoF}_h^{-1 / 2}\right)$ with different choice of $\kappa$. The associated adaptive mesh with $\kappa = 1$ is presented in  Fig. \ref{posterior_error and mesh example2} (b), which is generated after 8 refinement
iteration with 5,966 elements. It is clear to show that the algorithm does capture the solution feature and adaptively refines the mesh around the corners, where solution displays singularity.
In Table \ref{example2table}, numerical results are shown for the uniform and adaptive mesh refinements with  $\kappa=1$. It can be observed that the adaptive mesh refinement requires fewer $\mathrm{DoF}_h$ than the uniform mesh refinement to achieve greater accuracy.
\newline

\textbf{Example 3.} 
This example concerns the scattering by a grating surface with multiple sharp angles and also does not have an analytical solution. The parameters are set as $\mu=3, \lambda=2, \omega=1, \rho=1, \rho_f=1$ and   $\theta=\pi / 4 $. We take the period $\Lambda=5$.

\begin{figure}
	\centering
	\includegraphics[width=0.5\textwidth]{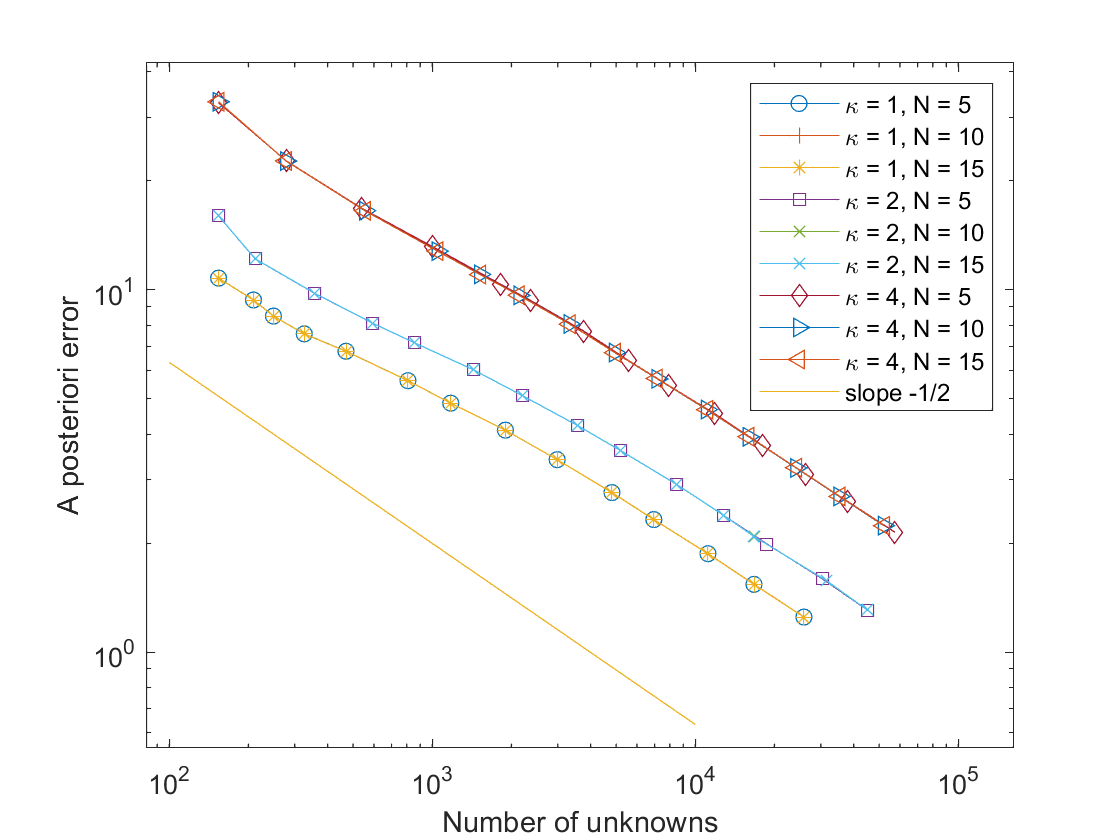}
	\caption{Quasi-optimality of the a posteriori errors with different with different truncated parameter $N$ in Example 3}
	\label{N}
\end{figure}

Fig \ref{k=1 example 3} shows the real parts of the numerical solutions with $\kappa = 1$.
Fig. \ref{posterior_error and mesh} (a) presents the curves of $\log \epsilon_h$ versus $\log \mathrm{DoF}_h$ with different wave number $\kappa$, where the decays of the a posteriori errors are $\mathcal{O}\left(\mathrm{DoF}_h^{-1 / 2}\right)$ with different choice of $\kappa$. The associated adaptive mesh with $\kappa = 1$ is presented in  Fig. \ref{posterior_error and mesh} (b), which is generated after 11 refinement
iteration with 11,168 elements. Once again,  the algorithm shows the ability to capture the singularity of the solution and perform additional local refinements. Fig. \ref{N} shows the curve of $\log \epsilon_h$ versus $\log \mathrm{DoF}_h$ with different truncated parameter $N$, which proves that
our adaptive algorithm is robust with respect to the choice of $N$. In Table \ref{example3table}, numerical results are shown for the uniform and adaptive mesh refinements with  $\kappa=1$, which shows the advantage of using adaptive mesh refinements.
\newline

\begin{figure}
	\centering
	\subfigure[A posterior error]{\includegraphics[width=0.495\textwidth]{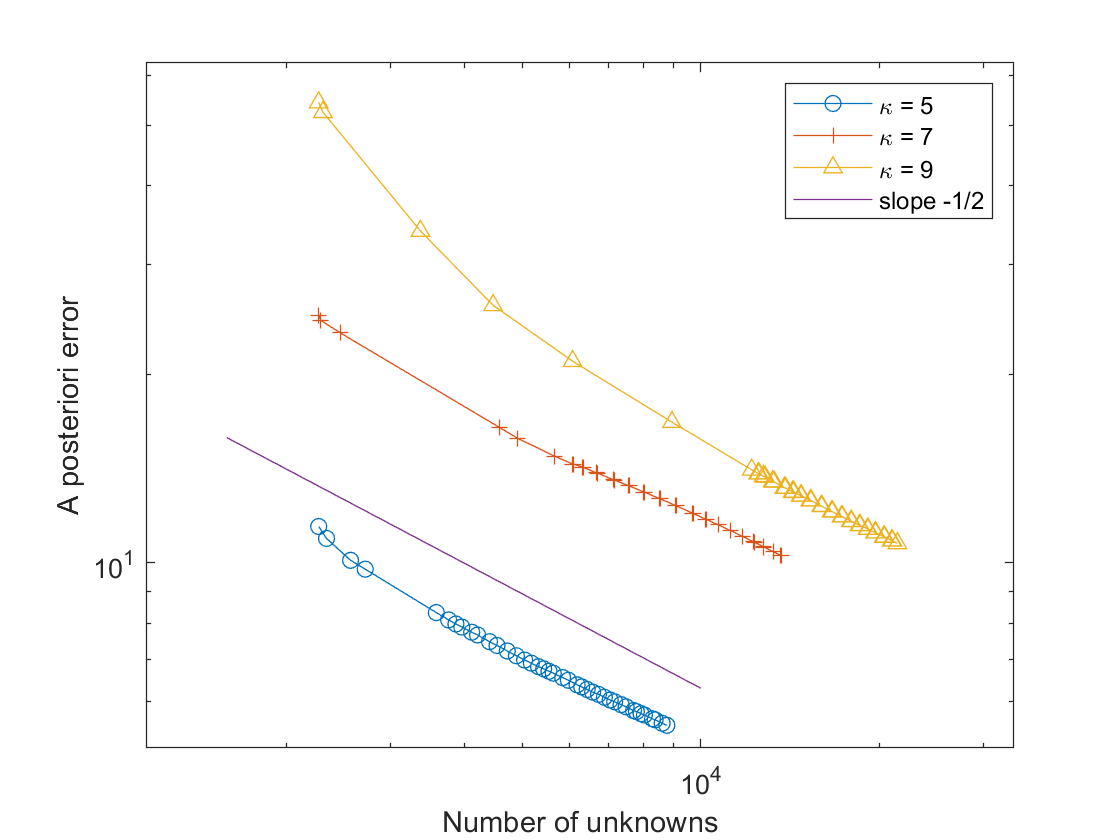}}
	\subfigure[Adaptive mesh]{\includegraphics[width=0.495\textwidth]{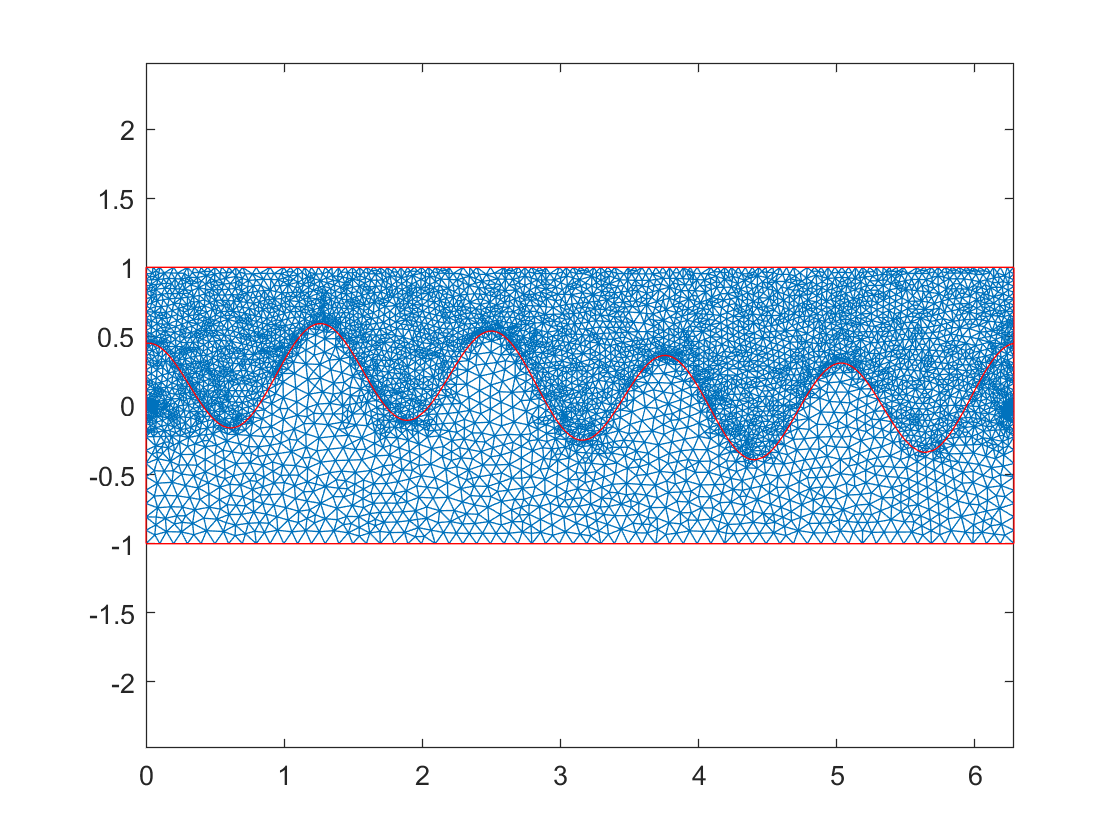}}
	\caption{Quasi-optimality of the a posteriori errors with different $\kappa$ (a); the adaptive mesh with 14,024 elements after 30 refinement iterations (b) in Example 4}
	\label{posterior_error and mesh exm4}
\end{figure}

\begin{table}
	\centering
	\caption{Comparison of numerical results using uniform and adaptive refinements with  $\kappa = 9$ in Example 4}
	\begin{tabular}{lllll}
		\hline \multicolumn{2}{l}{ Uniform mesh } && \multicolumn{2}{l}{Adaptive  mesh } \\
		\hline
		$\mathrm{DoF}_h$ & $\epsilon_h$  && $\mathrm{DoF}_h$  & $\epsilon_h$   \\
		\hline
		591 & 64.4435 && 591 & 64.4435 \\
		2265 & 54.3031 && 1014 & 49.2024  \\
		8889 & 21.6240 && 6077 & 20.9892  \\
		35217 & 10.7576 && 21485 & 10.7450  \\
		\hline
	\end{tabular}
	\label{example4table}
\end{table}	
	
\textbf{Example 4.} 
This example concerns the scattering by a curved grating surface. The parameters are set as $\mu=4, \lambda=2, \omega=1, \rho=1, \rho_f=1$ and   $\theta=\pi / 5 $. We take the period $\Lambda=2\pi$. The grating profile is given by
\begin{equation*}
	f(x_1) = 0.1+0.15\sin(x_1)+0.35\cos(5x_1),\quad x_1 \in [0,2\pi].
\end{equation*} 

Fig. \ref{posterior_error and mesh exm4} (a) presents the curves of $\log \epsilon_h$ versus $\log \mathrm{DoF}_h$, where the decays of the a posteriori errors are still quasi-optimal with different wave number $\kappa$. The associated adaptive mesh with $\kappa = 5$ is presented in  Fig. \ref{posterior_error and mesh exm4} (b), which is generated after 30 refinement
iteration with 14,024 elements. The algorithm shows the ability to capture the singularity around the fluid-solid interface and perform local refinements. In Table \ref{example4table}, numerical results are shown for the uniform and adaptive mesh refinements with  $\kappa=9$, which shows the advantage of using adaptive mesh refinements again.

\section{Conclusion}\label{conclusion}
In this paper, we present an adaptive finite element DtN method for the acoustic-elastic interaction
problem in periodic structures. A  duality argument
is applied to obtain the a posteriori error estimate. It does not only take into account
of the finite element discretization error but also includes the truncation error of two different wave
DtN operators. The theoretical analysis shows that the truncation error of the DtN
operators decays exponentially with respect to
the truncation parameter. The a posteriori error estimate for the discrete
problem serves as a basis for the adaptive finite element approximation. Numerical
results show that the proposed method is robust and effective. 
This work provides a viable alternative to the adaptive finite element
PML method for solving the acoustic-elastic interaction problem in periodic structures. It also enriches the range of choices available
for solving wave propagation problems imposed in unbounded domains.
Future work includes extending our analysis to the adaptive DtN-FEM for solving the electromagnetic-elastic interaction problem.

%
%
%

\bmhead{Acknowledgments}

This work of J.L. was partially supported by 
the National Natural Science Foundation of China grant 12271209.

\section*{Data availability statement}
All data generated or analysed during this study are included in this article.

\section*{Declarations}

\textbf{Conflict of interest} The authors declare no conflict of interest.

\end{document}